\documentclass[12pt,a4paper]{amsart}

\makeatletter
\long\def\unmarkedfootnote#1{{\long\def\@makefntext##1{##1}\footnotetext{#1}}}
\makeatother

\usepackage{amsmath,amssymb}%,a4}
\usepackage{enumerate,amssymb}
\usepackage{color}

\usepackage{graphicx}			% images and drawing
\usepackage{tikz}
\usetikzlibrary{cd}

\usepackage{hyperref}

\theoremstyle{plain}
\newtheorem{thm}{Theorem}[section]

\newtheorem{lemma}[thm]{Lemma}

\newtheorem{corollary}[thm]{Corollary}

\newtoks\prt
\newtheorem{proclaim}[thm]{\the\prt}
\theoremstyle{definition}

\newtheorem{remark}[thm]{Remark}

\newtheorem{definition}[thm]{Definition}

\def\eqn#1$$#2$${\begin{equation}\label#1#2\end{equation}}

\numberwithin{equation}{section}

\newcommand{\abs}[1]							% Absolute value
{\left| #1 \right|}

\def\diam{\operatorname{diam}}
\def\dist{\operatorname{dist}}

\def\epsilon{\varepsilon}
\def\en{\mathbb N}
\def\er{\mathbb R}

\def\haus{\mathcal{H}}

\def\loc{\operatorname{loc}}

\def\mir1{\mathcal L_1}

\def\oint{-\hskip -13pt \int}

\def\phi{\varphi}

\def\rn{\mathbb R^n}

\newtoks\by
\newtoks\paper
\newtoks\book
\newtoks\jour
\newtoks\yr
\newtoks\pages
\newtoks\vol
\newtoks\publ
\def\ota{{\hbox\vol{???}}}
\def\cLear{\by=\ota\paper=\ota\book=\ota\jour=\ota\yr=\ota
\pages=\ota\vol=\ota\publ=\ota}
\def\endpaper{\the\by, {\the\paper},
\textit{\the\jour} \textbf{\the\vol} (\the\yr), \the\pages.\cLear}
\def\endbook{\the\by, \textit{\the\book}, \the\publ.\cLear}
\def\endprep{\the\by, \textit{\the\paper}, \the\jour.\cLear}
\def\endyearprep{\the\by, \textit{\the\paper}, \the\jour, (\the\yr).\cLear}
\def\name#1#2{#2 #1}
\def\nom{ \rm no. }

\title{$(INV)$ condition and regularity of the inverse}

\author[A. Dole\v{z}alov\'a]{Anna Dole\v{z}alov\'a}
\address{Department of Mathematics and Statistics, University of
Jyv\"askyl\"a, P.O. Box 35 (MaD), 40014 Jyv\"askyl\"a, Finland
and
Department of Decision-Making Theory, Institute of Information Theory and Automation, Czech Academy of Sciences, Pod Vodárenskou věží 4, 182 00 Prague 8, Czech Republic}
\email{dolezalova@utia.cas.cz}

\author[S. Hencl]{Stanislav Hencl}
\address{Department of Mathematical Analysis, Faculty of Mathematics and Physics, Charles University, Sokolovsk\'a 83, 186 00 Prague 8, Czech Republic}
\email{hencl@karlin.mff.cuni.cz}

\author[J. Onninen]{Jani Onninen}
\address{Department of Mathematics, Syracuse University, Syracuse,
NY 13244, USA and  Department of Mathematics and Statistics, University of Jyv\"askyl\"a, P.O.Box 35 (MaD), 40014 Jyv\"askyl\"a, Finland
}
\email{jkonnine@syr.edu}

\thanks{A. Doležalová was supported by the Academy of Finland project 334014, and by the Czech Academy of Sciences project PPLZ L100752451. S. Hencl was supported by the grant GA\v CR P201/24-10505S. J. Onninen was supported by the NSF grant DMS-2453853.\\
Published at \url{doi.org/10.1016/j.jfa.2025.111215}.}

\date{\today}

\begin{document}

\begin{abstract}
Let $f \colon \Omega \to \Omega' $ be a Sobolev mapping of finite distortion between planar domains $\Omega $ and $\Omega'$, satisfying the $(INV)$ condition and coinciding with a homeomorphism near $\partial\Omega $. We show that $f$ admits a generalized inverse mapping $h \colon \Omega' \to \Omega$, which is also a Sobolev mapping of finite distortion and satisfies the $(INV)$ condition.  

We also establish a higher-dimensional analogue of this result: if  a mapping  $f \colon \Omega \to \Omega' $ of finite distortion is in the Sobolev class $W^{1,p}(\Omega, \mathbb{R}^n)$ with $p > n-1$ and satisfies the $(INV)$ condition, then $f$ has an inverse in $W^{1,1}(\Omega', \mathbb{R}^n)$ that is also of finite distortion.

Furthermore,  we characterize Sobolev mappings satisfying $(INV)$ whose generalized inverses have finite $n$-harmonic energy.  
\end{abstract}

\maketitle

\section{Introduction}
We study mappings fundamental to Geometric Function Theory (GFT) through variational approaches~\cite{AIMb, IMb, HK, Reb} which also serve as key models for deformations in mathematical models of Nonlinear Elasticity (NE). 
Let $\Omega\subseteq\rn$ be a domain that corresponds to some body that we want to deform and let $f:\Omega\to\rn$ be a mapping which describes our deformation of $\Omega$. Starting with the pioneering works of Ball \cite{Ba}, \v{S}ver\'ak \cite{Sv} and Ciarlet and Ne\v{c}as \cite{CN} we want our mappings to be injective or at least injective a.e. as it corresponds to the {\it non-interpenetration of matter}.

However, a serious challenge arises when passing to the weak limit of an energy-minimizing sequence of homeomorphisms $f_j \colon \Omega \to \Omega'$ in  $W^{1,p} (\Omega, \mathbb R^n)$: injectivity is generally lost. To develop viable mathematical models of hyperelasticity,  we find ourselves forced to adopt such limits as legitimate deformations in a way complying (as much as possible) with the principle of non-interpenetration of matter. 

When $p\ge n$, such weak limits (being also uniform limits) are monotone Sobolev mappings.  {\it Monotonicity} of a continuous mapping $f$, the concept of  Morrey~\cite{Mor}, simply means that  the preimage $f^{-1} (y_0)$ of any point  $ y_0 \in \mathbb R^n$ is a continuum. In the planar case, the weak $W^{1,p}$-limits of  homeomorphisms, when $p\ge 2$, are exactly monotone Sobolev mappings~\cite{IOmono}. This motivates us to adopt  monotone Sobolev mappings as legitimate deformations in the mathematical description of  plates. 

The case $p<n$, however, is significantly different: the limit of Sobolev homeomorphisms in $W^{1,p} (\Omega, \mathbb R^n)$ may fail to be continuous, allowing for phenomena such as cavitation. Consequently, monotone maps are no longer sufficient as admissible deformations. Nevertheless,  NE models typically exclude only mappings deemed physically inappropriate, still allowing  \emph{weak interpenetration of matter}; roughly speaking, squeezing of portion of the material can occur, but not folding, and no positive volume of material is permitted to turn "inside out." Building on an earlier work of Šverák~\cite{Sv}, this perspective led Müller and Spector~\cite{MS} to define the important concept of Sobolev mappings satisfying the so-called $(INV)$ condition. 

Informally, the $(INV)$ condition means that a ball $B(x,r)$ is mapped inside the image of its boundary $f(S(x,r))$, while the complement $\Omega\setminus \overline{B(x,r)}$ is mapped outside $f(S(x,r))$ (see Preliminaries for the formal definition). 
If $p>n-1$, then weak limits of  Sobolev homeomorphisms in $W^{1,p} (\Omega, \mathbb R^n)$  satisfy the $(INV)$ condition~\cite{MS}. Even more, the class of $(INV)$-class in $W^{1,p}$ is closed with respect to the weak topology when $p>n-1$, see~\cite[Lemma 3.3]{MS}. If  such deformations 
 have  strictly positive Jacobian a.e., then they are injective a.e.~\cite{MS}. For more recent results concerning the 
$(INV)$-class we refer the reader to e.g. \cite{BHMC, BHMCR2, CDL, DHM, DHMo, HMC, HeMo11, K, MST,SwaZie2002,SwaZie2004,T}. 

A longstanding open problem in NE is the existence of a homeomorphic minimizer for neohookean-type energy functionals, where the integrand blows up as the Jacobian approaches zero, see~\cite{Ba, BP,  BPO1, CN, CDL,   FG_95, MS, MST, Sv}.  In GFT, a related challenge is proving the existence of homeomorphic minimizers for the $L^p$-mean distortion energy when $p>1$,  i.e., to show that there always exists an $L^p$-Teichm\"uller mapping when $p>1$, see~\cite{IMO14, IMO, MY22, MY24}.  A natural approach is to construct a candidate for the inverse map and analyze its behavior by examining its associated inner-variational equation. However, determining whether such a generalized inverse even belongs to a Sobolev class is already a nontrivial issue. Indeed, even the inverse of a Sobolev homeomorphism does not necessarily belong to a Sobolev space.

 In the one-dimensional case, every homeomorphism and its inverse have bounded variation. A similar duality holds in the plane: a planar homeomorphism $f \colon \Omega \to \Omega'$ has bounded variation if and only if the inverse $f^{-1} \colon \Omega' \to \Omega$ has bounded variation, see~\cite{HKO}. Assuming that a homeomorphism $f \colon \Omega \to \Omega'$ has {\it finite distortion}  is the least regularity to imply that  the inverse $f^{-1}$  is  a Sobolev mapping.  Recall that a Sobolev homeomorphism  $f \in W_{\loc}^{1,1} (\Omega, \mathbb R^n)$ has  finite distortion if  there is a measurable function $K (x) \ge 1$, finite a.e., such that
\begin{equation}\label{eq_dist}
\abs{Df(x)}^n \le   K(x) \, J_f(x)\, , \qquad \text{where } J_f(x)=\det Df(x) .
\end{equation}
The smallest measurable $ K(x) \geq 1$ that satisfies the inequality~\eqref{eq_dist} is denoted by $K_f(x)$ and is referred to as the {\it distortion function} of  $f$.
A planar homeomorphism has finite distortion if and only if  its inverse has finite distortion~\cite{HK2}.

Our main result is an analogue of this foundational result for the $(INV)$ condition.

\prt{Theorem}
\begin{proclaim}\label{main}  
Let $\Omega $ and $\Omega' $ be domains in $ \mathbb{R}^2 $, and let $f_0 \colon \Omega \to \Omega' $ be a given homeomorphism. Suppose $f \colon \Omega \to \mathbb{R}^2$ is a mapping that coincides with $f_0$ on a neighborhood of $ \partial \Omega $. If  $f \in W^{1,1}_{loc}(\Omega, \mathbb{R}^2)$ has finite distortion and satisfies the $(INV)$ condition, then there exists a generalized inverse mapping $h \colon \Omega' \to \mathbb{R}^2$ with the same properties:  
\begin{itemize}
\item  $h \in W^{1,1}_{\loc}(\Omega', \mathbb{R}^2)$,  
\item  $h$ has finite distortion, and  
\item $ h$ satisfies the $(INV)$ condition.  
\end{itemize}
Furthermore,  the generalized inverse mapping  $h \colon \Omega' \to \Omega $ coincides with the classical inverse mapping a.e.\ on the set $f(\{x \in \Omega \colon J_f(x) > 0\})$. Additionally, the mappings satisfy  $x\in h_T(f(x))$ for a.e.\ $x\in\Omega$ and $y\in f_T(h(y))$ for a.e. $y\in\Omega'$.
\end{proclaim}  
Here $h_T (y)$ denotes the topological image of point $y\in \Omega'$ under the map $h$ (see~\eqref{imagepoint} for the precise definition). Note that in the planar case, a mapping $f \in W_{\loc}^{1,p}(\Omega, \mathbb{R}^2)$ satisfying $(INV)$ with $J_f > 0$ a.e. is invertible a.e. also in the limiting case $ p = 1$, and therefore the inverse mapping $f^{-1} \colon f(\{x \in \Omega \colon J_f(x) > 0\}) \to \{x \in \Omega \colon J_f(x) > 0\}$ is well-defined a.e. on the set.

The only previously known results related to Theorem~\ref{main} are due to Henao and Mora-Corral \cite[Theorem 5.1]{HMC} and \cite[Theorem 3.3]{HMC2}. However, their conclusions rely on significantly stronger assumptions, namely that $J_f > 0$  a.e. and that the Jacobian coincides with the distributional Jacobian, providing additional regularity properties such as the Luzin $(N)$ condition (for every set $E \subseteq \Omega $ with $|E| = 0$, we have $|f(E)| = 0$).  

In contrast, in our setting, null sets may be mapped to sets of positive measure, or null sets in the image may have preimages with positive measure (as the set $\{J_f = 0\}$ may have positive measure).  Nevertheless, the conclusion that the inverse mapping also satisfies the $(INV)$ condition is entirely new, even under the stronger assumptions considered by Henao and Mora-Corral~\cite{HMC, HMC2}. We believe that our main theorem is more suited to applications in GFT, where it is often unrealistic to assume $ J_f > 0$ a.e. or that the distributional Jacobian coincides with the pointwise Jacobian.   Furthermore, we provide a counterexample demonstrating that without the critical assumption of finite distortion, no reasonable theory of the $(INV)$-class can be developed (see Section~\ref{sec:counterex}). In such cases, nested balls may fail to have nested topological images, and the degree of a mapping can change arbitrarily (e.g., from $1$ to $-1$), undermining the fundamental structure necessary for this theory.

In higher dimensions, the inverse of a homeomorphism \( f \colon \Omega \to \Omega' \) of finite distortion, belonging to \( W_{\mathrm{loc}}^{1,n-1} (\Omega, \mathbb{R}^n) \), also has finite distortion. In particular, \( f^{-1} \in W_{\mathrm{loc}}^{1,1} (\Omega', \mathbb{R}^n) \) (see~\cite{CHM}). We establish the following analogous result.

\prt{Theorem}
\begin{proclaim}\label{main2}
Let $ \Omega $ and $\Omega' $ be domains in $ \mathbb{R}^n $, $n\geq 3$ and $p> n-1$, and let $f_0 \colon \Omega \to \Omega'$ be a given homeomorphism. Suppose that a mapping $ f \in W^{1,p}_{\mathrm{loc}}(\Omega, \mathbb{R}^n)$ has finite distortion, satisfies the $(INV)$ condition, and coincides with $f_0$ on a neighborhood of $\partial \Omega $. Then there exists a generalized inverse $h \colon \Omega' \to \mathbb{R}^n$ with finite distortion. In particular, $h \in W^{1,1}_{\mathrm{loc}}(\Omega', \mathbb{R}^n)$. 

Moreover, the generalized inverse mapping $h \colon \Omega' \to \mathbb{R}^n $ is differentiable a.e. in $\Omega' $ and coincides with the classical inverse mapping a.e.\ on the set $f(\{x \in \Omega \colon J_f(x) > 0\})$. Additionally, the mappings satisfy $x \in h_T(f(x))$ for a.e.\ \(x \in \Omega\) and $y \in f_T(h(y))$ for a.e. $ y \in \Omega'$.
\end{proclaim}

We conjecture that it suffices to take $p = n-1$ in Theorem~\ref{main2}. However, under this weaker assumption, one cannot conclude that $h$ is differentiable almost everywhere.

If a planar Sobolev homeomorphism $f \colon \Omega \to \Omega'$ has integrable distortion, its inverse mapping $h = f^{-1}$ has finite Dirichlet energy and the following identity holds (see~\cite{HKO05}):  
\begin{equation}\label{eq_identity_2d}
\int_{\Omega'} |Dh(y)|^2 \, dy = \int_{\Omega} K_f(x) \, dx.
\end{equation}

The key insight is that minimizing the $L^1$-distortion functional naturally transitions to studying the energy of the inverse mapping, which simplifies the analysis by resulting in a convex variational integral.

An analogue of this result for the $(INV)$ condition reads as follows.

\prt{Theorem}
\begin{proclaim}\label{thm_cor_1}  
Let $\Omega $ and $ \Omega' $ be bounded domains in $ \mathbb{R}^2 $, and let $f_0 \colon \Omega \to \Omega' $ be a given homeomorphism. Assume that $ f \in W^{1,1}_{\loc}(\Omega, \mathbb{R}^2)$ satisfies the $(INV)$ condition, coincides with $f_0$ on a neighborhood of $\partial \Omega $ and has integrable distortion. 

Then there exists a generalized inverse mapping \(h \colon \Omega' \to \mathbb{R}^2\) such that \(h \in W^{1,2}(\Omega', \mathbb{R}^2)\) and 
\[
\int_{\Omega'} |Dh(y)|^2 \, dy = \int_{\Omega} K_f(x) \chi_{G'}(x)\, dx.
\]
Here $G'= \{x\in \Omega \colon J_f(x)>0\}$. Moreover, $h$ satisfies the $(INV)$ condition.
\end{proclaim}
Note that any mapping $h \in W_{\loc}^{1,2}(\Omega', \mathbb{R}^2)$ satisfying the $(INV)$ condition admits a continuous representative. Moreover, if $h$ is a homeomorphism on a neighborhood of the boundary \(\partial \Omega'\), then it is monotone in the sense of Morrey~\cite[Lemma 2.8]{DPP}.

To establish that the inverse $h$ has finite $n$-harmonic energy (belongs to $W^{1,n}$) through a higher-dimensional analogue of the identity~\eqref{eq_identity_2d}, one must look at the inner distortion function of $f$. Indeed, it  can be shown easily (at least formally), that the pullback of the $n$-form $K_I(x) \, dx$ by the inverse mapping $h \colon \Omega' \to \Omega$ equals $|Dh(y)|^n \, dy \in \wedge^n \Omega'$.

For the precise formulation, recall that for a mapping $f \in W^{1,1}_{\mathrm{loc}}(\Omega, \mathbb{R}^n)$ of finite distortion, the {\it inner distortion function} $ K_I(x) = K_I(x, f) \geq 1 $ is defined as the smallest function satisfying  
$$
|D^\sharp f(x)|^n = K_I(x) \cdot J_f(x)^{n-1}.
$$
Here, the cofactor matrix  $D^\sharp f$ (also called the co-differential of $f$) encodes information about the $(n-1) \times (n-1)$-minors of the differential matrix  $Df$. Now, the precise result reads as follows (see ~\cite{AIMO, CHM, On06}):  

Let $ f \in W^{1,n-1}_{\mathrm{loc}}(\Omega, \mathbb{R}^n)$ be a homeomorphism of integrable inner distortion between bounded domains $\Omega, \Omega' \subseteq \mathbb{R}^n$. Then the inverse mapping $h = f^{-1}$ belongs to $W^{1,n}(\Omega', \mathbb{R}^n)$ and we have
$$
\int_{\Omega'} |Dh(y)|^n \, dy = \int_{\Omega} K_I(x, f) \, dx.
$$
We establish the following variant for mappings that satisfy only the $(INV)$ condition.
\prt{Theorem}
\begin{proclaim}\label{thm_cor_2}  
Let $\Omega $ and $ \Omega' $ be bounded domains in $ \mathbb{R}^n$, $n\geq 3$ and $p>n-1$, and let $f_0 \colon \Omega \to \Omega' $ be a given homeomorphism. Assume that $ f \in W^{1,p}_{\loc}(\Omega, \mathbb{R}^n)$ satisfies the $(INV)$ condition, coincides with $f_0$ on a neighborhood of $\partial \Omega $ and has integrable inner distortion. 

Then there exists a generalized inverse mapping $h \colon \Omega' \to \mathbb{R}^n$ with  $h \in W^{1,n}(\Omega', \mathbb{R}^n)$ and
$$
\int_{\Omega'} |Dh(y)|^n \, dy = \int_{\Omega} K_I(x, f) \chi_{G'}(x) \, dx.
$$
Here $G'= \{x\in \Omega \colon J_f(x)>0\}$. 
\end{proclaim}

Let us briefly outline the structure of the paper. 

Following the {\bf Preliminaries}, we analyze in {\bf Section~\ref{sec:3}}  the $(INV)$ condition for mappings of finite distortion without the usual assumption that $J_f > 0$ a.e. Even without this key assumption, we ``reprove'' many standard conclusions: the mapping $f$ is injective a.e. on the set $G'$; it satisfies a stronger version of the $(INV)$ condition (all points of a ball are mapped inside the topological image, not just almost every point); it is differentiable almost everywhere; and so on. This framework also allows us to define the topological image of a point in the usual way (see \eqref{imagepoint} below) and to establish that there are at most countably many cavities, i.e., points $x\in\Omega$, for which the topological image  $f_T(x)$ has positive measure. 

In {\bf Section~\ref{sec:h_def_ae}}, we leverage these results to define the mapping $h$ on $\Omega'$. For cavities formed at a point $x$, we define $h(y) = x$ for $y \in f_T(x)$. For ``good'' points  $x$, where $f$ is injective, we naturally set $h(f(x)) = x$. We show that for each ball $B(y, r) \subseteq f(\Omega)$, $h$ is defined on a set of positive measure. This allows us to extend $h(y)$ to all points $y$ as a limit of integral averages over regions where $h$ is already defined. 

In {\bf Section~\ref{sec:reg_of_inv}}, we prove Theorems \ref{main}  (modulo that $h$ satisfies the $(INV)$ condition) and \ref{main2} for the mapping $h$ defined in the previous section.

In {\bf Section~\ref{sec:2d_h_inv}}, we show that in the planar case the generalized inverse $h$ also satisfies $(INV)$.

Finally, in {\bf Section~\ref{sec:conf_inner}}, we prove the connection between the integrable inner distortion and the finite  $n$-harmonic energy of its generalized inverse map 

To avoid confusion, for the rest of the paper, $f^{-1}(y)$ refers to the set of preimages of $y$ under $f$, whereas $h$ represents the generalized inverse mapping. While $h$ and $f^{-1}$ can be identified for mappings where $J_f > 0$ a.e., they may differ in the more general setting of mappings of finite distortion.

\section{Preliminaries}

\subsection{$(INV)$ condition}

Suppose that
$f\colon S(a,r) \to \mathbb{R}^n$ is continuous,
following \cite{MS} we define a {\it topological image} of $B(a,r)$ as
$$
im_T (f,B):=\bigl\{y\in \rn\setminus f(S(a,r)):\ \deg(f,S(a,r),y)\neq 0\bigr\}.
$$
Denote
$$E(f,B(a,r)):=im_T (f,B(a,r))\cup f(S(a,r)).$$

\begin{definition}
	Let $f:\Omega\to\mathbb{R}^n$. We say that $f$ satisfies $(INV)$ for a ball 
	$B\subset\subset\Omega$ 
	if 
	\begin{enumerate}
		\item it is continuous on $\partial B$;
		\item $f(x)\in im_T (f,B)\cup f(\partial B)$ for a.e.\ $x\in \overline{B}$;
		\item $f(x)\notin im_T (f,B)$ for a.e.\ $x\in\Omega\setminus B$.
	\end{enumerate}
	We say that $f$ satisfies the $(INV)$ condition if for every $a\in\Omega$ there is $r_a>0$ such that for $\mathcal{H}^1$-a.e.~$r\in (0,r_a)$ it satisfies $(INV)$ for $B(a,r)$. 
\end{definition}

Let us recall that a quasicontinuous representative of $f\in W^{1,p}(\Omega,\rn)$, $p>n-1$, is continuous for every $x$ on almost every sphere $S(x,r)$ (see e.g. \cite[Theorem 3.3.3 and Theorem 2.6.16]{Z}).

\subsection{Density points}
\begin{definition}
	Let $A\subseteq \mathbb{R}^n$ be a measurable set and $x\in\mathbb{R}^n$. We denote by 
	$$
	D^+(x,A)=\limsup_{r\to 0+} \frac{|B(x,r)\cap A|}{|B(x,r)|}
	$$ 
	and
	$$
	D^-(x,A)=\liminf_{r\to 0+} \frac{|B(x,r)\cap A|}{|B(x,r)|}
	$$ 
	the upper and lower density of $A$ around $x$. If those two values coincide, we denote their value by $D(x,A)$ and call it the density of $A$ around $x$. We can also say for $D(x,A)=1$ that $x$ is a point of density (one) of $A$. 
\end{definition}
It is well-known that a.e.~$x\in A$ is a point of density of $A$. 

\subsection{Change of variables}

The following statements could be found e.g. in \cite[Theorem A.35]{HK}. 

\prt{Theorem}
\begin{proclaim}\label{subst}
Let $f\in W_{\loc}^{1,1}(\Omega,\rn)$ and let $\eta$ be a nonnegative Borel measurable function on $\rn$. 
Then 
\eqn{sub}
$$
\int_{\Omega}\eta(f(x))\,|J_f(x)|\,dx\le \int_{\rn}\eta(y)\,N(f,\Omega,y)\,dy,
$$
where the multiplicity function $N(f,\Omega,y)$ of $f$ is defined as the number of preimages of
$y$ under $f$ in $\Omega$. Moreover, there is an equality in \eqref{sub} if we assume in addition that $f$ satisfies the Luzin $(N)$ condition. 
\end{proclaim}

\subsection{ACL condition}

The following statement could be found e.g. in \cite[Lemma A.15]{HK}. 

\prt{Theorem}
\begin{proclaim}\label{ACL}
Let $\Omega\subseteq\rn$ be a domain and let $u\in W^{1,p}(\Omega)$. Then $u$ has a representative $\tilde{u}$ that is absolutely continuous on almost all line segments in $\Omega$ parallel to the coordinate axes and whose (classical) partial derivatives belong to $L^p(\Omega)$.
\end{proclaim}

\subsection{Embedding into H\"older spaces} 

Next, we recall a well-known fact that a
function in the Sobolev class $W_{\loc}^{1,p} (\Omega)$, where $\Omega \subseteq \rn $, is H\"older continuous
with exponent $1-n/p$ provided $p>n$. More precisely, we have the following
oscillation lemma (see for example~\cite[Lemma 2.19]{HK}).

\begin{lemma}\label{lem:soboembed}
	Let $u$ be a function belonging to the Sobolev class $W_{\loc}^{1,p} (\Omega)$ where $\Omega \subseteq \rn $ and $p>n$. Then
	\begin{equation}\label{eq:soboembed}
		\abs{u(x) - u(y)} \le C(p,n) s^{1-\frac{n}{p}} \left(\int_{B_s}  \abs{D u}^p\right)^\frac{1}{p}
	\end{equation}
	for almost every $x,y \in B_s=B(x,s) \subseteq \Omega$. If $n=1$, then the estimate~\eqref{eq:soboembed} also holds for $p=n=1$. Furthermore, for the continuous representative the estimate~\eqref{eq:soboembed} holds for every $x,y \in B_s$. 
\end{lemma}

\subsection{Hausdorff measure and $p$-capacity}

The following statements follow e.g. from \cite[Chapter 4.7.2]{EG} (see also \cite[Theorem B]{A} for the case $n=2$ and $p=1$).

\begin{lemma}\label{hausdorff_implies_capacity}
Let $1\leq p < n$ and $P\subseteq \mathbb{R}^n$. Then 
$$
\mathcal{H}^{n-p}(P) = 0 \implies cap_p(P)=0.
$$
\end{lemma}

We also need the other implication.
\begin{lemma}\label{capacity_implies_hausdorff}
Let $1\leq p <n$ and $P\subseteq \mathbb{R}^n$. Then 
$$
cap_p(P)=0 \implies\mathcal{H}^{s}(P) = 0
$$
for all $s>n-p$.  Moreover, if $n=2$ and $p=1$, then it holds also for $s=1$.
\end{lemma}

\section{The mappings satisfy the strong $(INV+)$ condition}\label{sec:3}

\subsection{Variations of the $(INV)$ condition}

We introduce here different modifications of the $(INV)$ condition. Our ultimate goal is to show that for $f\in W^{1,p}_{loc}(\Omega, \mathbb{R}^n)$, $p>n-1$ (or $p\geq n-1$ for $n=2$), of finite distortion and satisfying $(INV)$ we can find a representative which satisfies the strong $(INV+)$ condition.

\begin{definition}
	Let $f:\Omega\to\mathbb{R}^n$ and $p\in [1,n)$. We say that $f$ satisfies the $(INV)_p$ condition if there exists a set $P\subseteq \Omega$ of $p$-capacity zero such that for every $a\in\Omega$ there is $r_a>0$ such that for $\mathcal{H}^1$-a.e.~$r\in (0,r_a)$ we have 
	\begin{enumerate}
		\item $f$ is continuous on $\partial B(a,r)$;
		\item $f(x)\in im_T (f,B(a,r))\cup f(\partial B(a,r))$ for every $x\in \overline{B(a,r)}\setminus P$;
		\item $f(x)\notin im_T (f,B(a,r))$ for every $x\in\Omega\setminus (B(a,r)\cup P)$.
	\end{enumerate}
\end{definition}

\begin{definition}
	Let $f:\Omega\to\mathbb{R}^n$. We say that $f$ satisfies strong $(INV)$ for a ball 
	$B\subset\subset\Omega$ 
	if 
	\begin{enumerate}
		\item $f$ is continuous on $\partial B$;
		\item $f(x)\in im_T (f,B)\cup f(\partial B)$ for every $x\in \overline{B}$;
		\item $f(x)\notin im_T (f,B)$ for every $x\in\Omega\setminus B$.
	\end{enumerate}
	We say that $f$ satisfies the strong $(INV)$ condition if for every $a\in\Omega$ there is $r_a>0$ such that for $\mathcal{H}^1$-a.e.~$r\in (0,r_a)$ it satisfies strong $(INV)$ for $B(a,r)$. 
\end{definition}

\begin{definition}
	We say that $f:\Omega\to\mathbb{R}^n$ satisfies the $(INV+)$ (or $(INV+)_p$, or strong $(INV+)$) condition, if $f\circ \Phi$ satisfies the $(INV)$ (or $(INV)_p$, or strong $(INV)$) condition for every bi-Lipschitz map $\Phi:\Omega \to\Omega$.
\end{definition}

\begin{lemma}\label{lemma_p_INV}
Let $f\in W^{1,p}_{loc}(\Omega, \mathbb{R}^n)$, $1\leq p<n$, such that $f^*$ satisfies $(INV)$, resp. $(INV+)$. Then $f^*$ satisfies $(INV)_p$, resp. $(INV+)_p$.
\end{lemma}
Remember that the precise representative $f^*(x)=\lim_{r\to 0+} \oint_{B(x,r)}f(z)dz$ whenever the limit exists (and $0$ otherwise).

\begin{proof}
The proof for $(INV)_p$ can be found in \cite[Lemma 7.2]{MS} and is not dependent on whether the set in question is a ball or a Lipschitz domain, as the important sets are closed anyway and the precise representative is approximatively continuous on $\Omega$ up to a set of $p$-capacity zero.
\end{proof}

Next we construct a representative which satisfies not only $(INV)$, but also strong $(INV)$.

\begin{lemma}\label{thm_injectivity}
Let $\Omega\subseteq \mathbb{R}^n$ be a domain and $f\in W^{1,p}_{loc}(\Omega, \mathbb{R}^n)$, $p>n-1$ (or $p\geq n-1$ for $n=2$), be a mapping satisfying $(INV)$. Denote $E\subseteq \Omega$ the set where $f$ is approximatively differentiable and $J_f>0$. Then $f$ is injective a.e. on $E$.
\end{lemma}

\begin{proof}
We fix the representative of $f$. From \cite[Lemma 2.5]{MS} applied on set $E$ we get a set $E'\subseteq E$ with $|E|=|E'|$ such that for every $x\in E'$ and $A\subseteq E$ Borel we have
\begin{equation}\label{density}
D(x,A)=1 \implies D(f(x), f(A))=1.
\end{equation}
Now we mimic the proof of \cite[Proposition 3.4]{MS}.
For contradiction, assume that we have two distinct points $x_1, x_2 \in E'$ such that $f(x_1)=f(x_2)$. We can assume that the density of $E'$ around those points is 1 (that holds for a.e. point of a set). We can find small spheres $S_1=S(x_1,r)$ and $S_2=S(x_2,r)$ around these points such that $f$ is in $W^{1,p}$ on $S_1$ and $S_2$ and the encircled balls $B_1$, $B_2$ are disjoint.

Applying $(INV)$ gives us sets $N_1$ and $N_2$ of zero measure such that 
$$
f(B_1\cap E'\setminus N_1)\subseteq im_T(f,B_1) \cup f(S_1)
$$
 and 
$$
 f(B_2\cap E'\setminus N_2)\subseteq \mathbb{R}^n\setminus im_T(f,B_1).
 $$
 
Therefore, 
$$f(B_1\cap E'\setminus N_1)\cap f(B_2\cap E'\setminus N_2)\subseteq f(S_1).$$ 
By \cite[Theorem 1]{MM}, $\mathcal{L}^n(f(S_1))=0$. (In case $p=1$ and $n=2$ we have that since $f$ is absolutely continuous on the circle by the analogy of the ACL condition.)

However, we know from \eqref{density} that both $f(E'\cap B_1\setminus N_1)$ and $f(E'\cap B_2\setminus N_2)$ have density 1 around $f(x_1)=f(x_2)$, and so has their intersection. Thus we obtain a contradiction.
\end{proof}

In the following theorem we show that topological images of nested balls are nested and topological images of disjoint balls are disjoint. Similar statement could be shown also e.g. for images of cubes (and not balls) so we could use it later also for other shapes. 

\begin{thm}\label{thm_embedded_ball}
Let $\Omega\subseteq \mathbb{R}^n$ be a domain and $f\in W^{1,p}_{loc}(\Omega, \mathbb{R}^n)$, $p>n-1$ (or $p\geq n-1$ for $n=2$), be a mapping of finite distortion satisfying $(INV)$. We have that for every $x_1, x_2\in \Omega$ and almost all radii small enough
\begin{itemize}
\item $B(x_1,r_1)\subseteq B(x_2, r_2) \implies E(f,B(x_1,r_1)) \subseteq E(f,B(x_2, r_2)),$
\item $B(x_1,r_1)\cap B(x_2, r_2) =\emptyset \implies im_T(f,B(x_1,r_1))\cap im_T(f,B(x_2, r_2))=\emptyset$.
\end{itemize}
\end{thm}

\begin{proof}
We mimic the proof of \cite[Lemma 7.3]{MS}, taking care of the set $\{J_f=0\}$ separately. 

From Lemma \ref{lemma_p_INV} we have $(INV)_p$, the representative is fixed from now on.
Denote $P$ the set where $(INV)$ is broken; by Lemma \ref{capacity_implies_hausdorff} we know that $\mathcal{H}^{1}(P)=0$.

Now let's take $E'$ from Lemma \ref{thm_injectivity}, so $f$ is one-to-one on $E'$ and \eqref{density} holds for every $x \in E'$. 

For radii small enough, up to a $\mathcal{H}^1$-measure zero set $N_i$, we have 
\begin{itemize}
\item $f$ is in $W^{1,p}(S(x_i,r_i), \mathbb{R}^n)$ and continuous on the sphere,
\item $(INV)_p$ is satisfied on $S(x_i,r_i)$,
\item $|\Lambda_{n-1}(Df)\nu|>0$ $\mathcal{H}^{n-1}$-a.e. on $S(x_i,r_i)\cap E'$ ($\nu$ is the outer normal to $S(x_i,r_i)$),
\item $\mathcal{H}^{n-1}(S(x_i,r_i)\cap (E\setminus E'))=0$,
\item $S(x_i,r_i)\cap P=\emptyset$ and
\item $|Df|=0$ $\mathcal{H}^{n-1}$-a.e. on $S(x_i,r_i)\setminus E'$, since $f$ is a mapping of finite distortion.
\end{itemize}
Here $\Lambda_{n-1}(Df)$ denotes the restriction of $Df$ to the tangent hyperplane to $S(x_i,r_i)$ as usual.
Let's fix such $r_1$ and $r_2$ and denote $B_i=B(x_i, r_i)$ and $S_i=S(x_i, r_i)$.

Since $f$ is continuous on $S_i$ and $|f(S_i)|=0$, the classical degree and its generalized version from \cite{CDL} coincide. Therefore from \cite[Proposition 3.4]{CDL} (or \cite[Lemma 2.3]{DHM} for general $n$) it is a BV function which is constant on components of $\mathbb{R}^n \setminus f(S_i)$ - note that the assumption $J_f>0$ a.e. is not used anywhere there. The components of $\mathbb{R}^n \setminus f(S_i)$ are open sets with boundary contained in $f(S_i)$ and the degree is zero on the unbounded component. So
$A_i:=im_T (f,B_i)$ are bounded open sets of finite perimeter with $\partial A_i\subseteq f(S_i)$ and $$
C:=E(f,B_2)=A_2\cup f(S_2)
$$ 
is a compact set of finite perimeter. (We need $C$ to be a closed set so that its complement is open.)

Outside of the set $E'$ the derivative is $0$ $\mathcal{H}^{n-1}$-a.e., so by $p>n-1$ and \cite[Theorem 1]{MM} (or by absolutely continuity for $n=2, p=1$) we have that
\eqn{ttt}
$$
\mathcal{H}^{n-1}(f(S_i\setminus E')) = 0.
$$
As $f$ is injective on $E'$, we have 
\eqn{tttt}
$$
f(S_1\cap E')\cap  f(S_2\cap E') =\emptyset.
$$

Now let's prove $(i)$: Assume $B_1\subseteq B_2$. Since $S_1$ does not intersect $P$ and we have $(INV)_p$ for $B_2$, we obtain $f(S_1)\subseteq C$. Similarly by switching $S_1$ and $S_2$ we get $f(S_2)\subseteq\mathbb{R}^n\setminus A_1$.

Now $A_1$ is an open bounded domain, $\mathbb{R}^n\setminus C$ is also open, both of finite perimeter. Since 
$$
\partial A_1\subseteq f(S_1)\subseteq C \text{ and } \partial C\subseteq f(S_2)\subseteq\mathbb{R}^n\setminus A_1,
$$
we can use \eqref{ttt} and \eqref{tttt} to obtain 
$$
\haus^{n-1}\bigl(\partial A_1\cap \partial (\rn\setminus C)\bigr)\leq 
\haus^{n-1}(f(S_1\cap E')\cap  f(S_2\cap E'))+\haus^{n-1}(f(S_1\setminus E'))+\haus^{n-1}(f(S_2\setminus E'))=0.
$$
Now we can use \cite[Lemma A.1]{MS} for the sets $A_1$ and $\rn\setminus C$ to obtain $A_1\cap (\mathbb{R}^n\setminus C)=\emptyset$, i.e., $A_1\subseteq C$.

To prove $(ii)$, we do basically the same: Assume $B_1\cap B_2 = \emptyset$. Since we are not intersecting $P$ and have $(INV)_p$, we know that 
$$f(S_1)\subseteq\mathbb{R}^n\setminus A_2 \text{ and }f(S_2)\subseteq\mathbb{R}^n\setminus A_1.$$
We again use Lemma A.1 (for sets $A_i$) obtaining that $A_1\cap A_2=\emptyset$.
 
\end{proof}

Under the conditions of the previous statement, we define 
the topological image of point $x$ as
\eqn{imagepoint}
$$
f_T(x)=\bigcap_{r>0,r\notin N_x} E(f^*,B(x,r)).
$$

\begin{thm}
Let $\Omega\subseteq \mathbb{R}^n$ be a domain and $f\in W^{1,p}_{loc}(\Omega, \mathbb{R}^n)$, $p>n-1$ (or $p\geq n-1$ for $n=2$), be a mapping of finite distortion satisfying $(INV)$. Then there exist a set $NC$ of zero $p$-capacity and $\hat{f}$ a representative of $f$ which is continuous on $\Omega\setminus NC$, $\hat{f}=f^*$ cap$_p$-a.e. and $\hat{f}(x)\in f_T(x)$ for any $x\in \Omega$.
\end{thm}

\begin{proof}
The proof is identical to the one of \cite[Theorem 7.4]{MS}, as the assumption $J_f\neq 0$ is used only to get the result analogous to Theorem \ref{thm_embedded_ball}, and the formulas hold also for $p=1$, $n=2$.
\end{proof}

\begin{lemma}\label{lemma_strong}
Under the assumptions of the previous theorem, if $f$ satisfies $(INV)$, resp. $(INV+)$, then it has a representative which satisfies strong $(INV)$, resp. strong $(INV+)$.
\end{lemma}

\begin{proof}
The proof is identical to the one of \cite[Corollary 7.5]{MS}.
\end{proof}

\begin{lemma}\label{diff2}
Let $\Omega\subseteq \mathbb{R}^n$ be a domain and $f\in W^{1,p}_{loc}(\Omega, \mathbb{R}^n)$, $p>n-1$ (or $p\geq n-1$ for $n=2$), satisfy the strong $(INV)$ condition. Then $f$ is differentiable a.e.
\end{lemma}

\begin{proof}
We take $x\in \Omega$ a Lebesgue point of $|Df|^p$ and $r\in (0,r_x/2)$. For a.e. $t\in [r,2r]$ the strong $(INV)$ condition holds on $B(x,t)$. That gives us for such $t$ (by the Sobolev embedding theorem on spheres or by absolutely continuity in the case $n=2$, $p=1$)
$$
\operatorname{osc}_{B(x,r)}f \leq \operatorname{osc}_{B(x,t)}f\leq \operatorname{osc}_{S(x,t)}f \leq C t\left(t^{-n+1}\int_{S(x,t)} |Df|^p\right)^{1/p}.	
$$
Integrating over $(r,2r)$ yields
$$
\frac{\operatorname{osc}_{B(x,r)}f}{r} \leq  C \left(\oint_{B(x,2r)} |Df|^p\right)^{1/p}< \infty.	
$$
By the Stepanov theorem (see e.g. \cite[Theorem 2.23]{HK}), $f$ is differentiable in $x$.
\end{proof}

To conclude this subsection, we show that the $(INV)$ condition implies the $(INV+)$ condition also under the assumption of $f$ being a mapping of finite distortion and not only when we ask for $J_f>0$ a.e. Note that if we drop the condition on the distortion, the result is not true anymore, see \cite[Section 5]{DPP}.

\begin{lemma}\label{lem_degree}
Let $\Omega\subseteq \mathbb{R}^n$ be a domain and $f\in W^{1,p}_{loc}(\Omega, \mathbb{R}^n)$, $p>n-1$ (or $p\geq n-1$ for $n=2$), be a mapping of finite distortion satisfying $(INV)$. Let $a\in\Omega$. Then there exists $r_a>0$ such that for $\haus^1$-a.e. radius $r\in (0,r_a)$ we have 
$$
\deg(f, S(a,r), y)\in \{0,1\} \text{ for all } y\in\mathbb{R}^n\setminus f(S(a,r)).
$$
\end{lemma}

\begin{proof}
We note here only the key differences to the proof of \cite[Lemma 3.5(ii)]{MS}, as most of the proof runs the same way. Pick $r$ such that $f$ satisfies $(INV)$ for $B(a,r)$. If $J_f=0$ $\haus^{n-1}$-a.e. on the sphere $S=S(a,r)$ then we can assume that $|Df|=0$ $\haus^{n-1}$-a.e. on $S$ by finite distortion. Then have that its image is a point and the statement obviously holds. Let's assume that $J_f>0$ on a subset of $S$ of positive $\haus^{n-1}$ measure from now on.

In Step 1, we choose $r$ to have the same properties as had $r_1$ and $r_2$ in the proof of Theorem \ref{thm_embedded_ball} i.e., intersecting the sphere with the good set $E'$. In particular, we have $|Df|=0$ $\haus^{n-1}$-a.e. on the set where $J_f=0$. 
In Step 2, we set $\tilde{\nu}$ to be 0 outside of $f(E')$ (we know that it is a set of $\haus^{n-1}$-measure 0). The key ingredient of this step is the following identity at the end of Step 2 (see notation from \cite[Lemma 3.5(ii)]{MS})
$$
\int_S (g\circ f)\cdot(\Lambda_{n-1} Df)\nu \; d\haus^{n-1}=\int_{f(S)}g\cdot \tilde{\nu}\; d\haus^{n-1}.
$$
This identity holds under our assumptions as we are either on a set where $J_f>0$ and our $f$ is $1-1$ a.e. there and we can use the change of variables formula (Theorem \ref{subst}). Or $J_f=0$ and then $\Lambda_{n-1} Df=0$ a.e. (which gives zero on the left-hand side) and $\haus^{n-1}$ of the image of this set with $J_f=0$ (and thus $|Df|=0$ a.e.) is zero, which gives zero also on the right-hand side. 
In Step 4, we choose $V\subseteq S\cap E'$. The rest of the argument in \cite[Lemma 3.5(ii)]{MS} follows.
\end{proof}

\prt{Theorem}
\begin{proclaim}\label{thm_inv+}
Let $\Omega\subseteq\mathbb{R}^n$ be a domain and $f:\Omega\to\mathbb{R}^n$ be a mapping of finite distortion in $W^{1,p}$, $p>n-1$ (or $p\geq n-1 $ for $n=2$), that satisfies $(INV)$. Then $f$ satisfies $(INV+)$.
\end{proclaim}

\begin{proof}

We can mimic the proof of \cite[Theorem 9.1]{MS} using Lemma \ref{lem_degree} to prove the theorem for a diffeomorphism $\Phi$ up to the very last step of the proof (the extension to bi-Lipschitz mappings follows as in \cite{DPP}). There they obtain sets $M_1$ and $M_2$ in the target, which are of zero measure, and their preimages are the only points where $(INV+)$ can be broken. We want to show that the preimages are of measure zero also in our case.

Let's assume for contradiction that $N_1:=f^{-1}(M_1)$ (the subset of $U$ which is mapped outside of $E(f,U)$) has positive measure. Therefore $J_f=0$ a.e. on $N_1$ (by change of variables as $|M_1|=0$), and so $|Df|=0$ a.e. on $N_1$ since $f$ is of finite distortion. We can find a lot of segments $L$ (in the same direction) connecting $\partial U$ and $N_1$ such that $f|L$ is absolutely continuous and of finite distortion. Such segment $L$ is mapped to $f(L)$, which has to connect $E(f,U)$ with some point from $M_1$, i.e., it has to cross a non-trivial distance outside of $E(f,U)$. Therefore 
$$
\haus^1(f(L)\setminus E(f,U))>0.
$$
However, only points from $L\cap N_1$ can be mapped there, and $f(L\cap N_1)$ is of zero $\haus^1$-measure by the one-dimensional change of variables as $|Df|=0$ $\haus^1$-a.e. on $L\cap N_1$. The proof for $M_2$ is analogous. 
 
What remains is to show \cite[Lemma 8.1]{MS} in our setting. The proof until \cite[(8.7)]{MS} is the same. 
For $a\in G'$ (that is, $J_f(a)>0$) we can assume that $f$ is differentiable at $a$ by Lemma \ref{diff2} and proceed the same way as their proof goes. If $a\notin G$, we can assume that $|Df(a)|=0$ and $f$ is differentiable in $a$. Therefore from the definition of the derivative we obtain
$$
\lim_{r\to 0+} r^{-n}|im_T(f,B(a,r))|\leq C \varepsilon^n
$$
for an arbitrary $\varepsilon>0$, giving us the required analogy of their (8.10) also a.e. outside of $G$ and we can conclude.

Note that we are also using Lemma \ref{lem_degree} which gives us the same result as  \cite[Lemma 3.5(ii)]{MS}.

\end{proof}

\subsection{Counterexamples}\label{sec:counterex}

The following example shows that the condition that we have a mapping of finite distortion which satisfies the $(INV)$ condition is indeed crucial for our conclusions. Without it we do not know that nested balls have nested topological images (i.e., Theorem \ref{thm_embedded_ball} does not hold) or that the degree cannot become $-1$, and so there is no reasonable theory. For simplicity we construct the mapping to satisfy the $(INV)$ condition with respect to cubes and not balls but the difference is not important. 

\prt{Lemma}
\begin{proclaim}\label{noMFDlemma2}
There is a mapping $f:[-2,2]^2\to[-2,2]^2$ which is identity on the boundary, it belongs to $W^{1,p}([-2,2]^2,[-2,2]^2)$ for any $p<2$, it satisfies the $(INV)$ condition (with respect to cubes) but for cubes $Q=[-1,1]^2$ and $Q_1=[-\tfrac{1}{2},0]\times[\tfrac{1}{2},1]$ we have 
$$
Q_1\subseteq Q, f(Q_1)\subseteq f(\partial Q)\text{ but }im_T (f,Q_1)\text{ is not a subset of }im_T (f,Q).  
$$ 
Moreover,
$$
\deg(f,Q_1,y)=-1\text{ for every }y\in im_T (f,Q_1)
$$
and the same conclusion holds for any small $t>0$ and cubes $[-1-t,1+t]^2$ and $[-\tfrac{1}{2}+t,0-t]\times[\tfrac{1}{2}+t,1-t]$. 
\end{proclaim}
\begin{proof}
Let us denote (see Fig \ref{noMFD})
$$
Q=[-1,1]^2,\ Q_1=[-\tfrac{1}{2},0]\times[\tfrac{1}{2},1],\ Q_2=[0,\tfrac{1}{2}]\times[\tfrac{1}{2},1]\text{ and }
Q_3=[-\tfrac{1}{2},0]\times[1,\tfrac{3}{2}]. 
$$
We construct our mapping as a composition of several mappings 
$$
f=r\circ l\circ b\circ k
$$
and the overall idea is depicted in Fig. \ref{noMFD}. All our mappings are identity on $\partial [-2,2]^2$ so their composition $f$ satisfies this as well. Moreover, we construct them so that $b$ is bi-Lipschitz and $l$ and $r$ are Lipschitz. 

\begin{figure}[h t p]
	\vskip 210pt
	{\begin{picture}(0.0,0.0) 
			\put(-150.2,0.2){\includegraphics[width=0.90\textwidth]{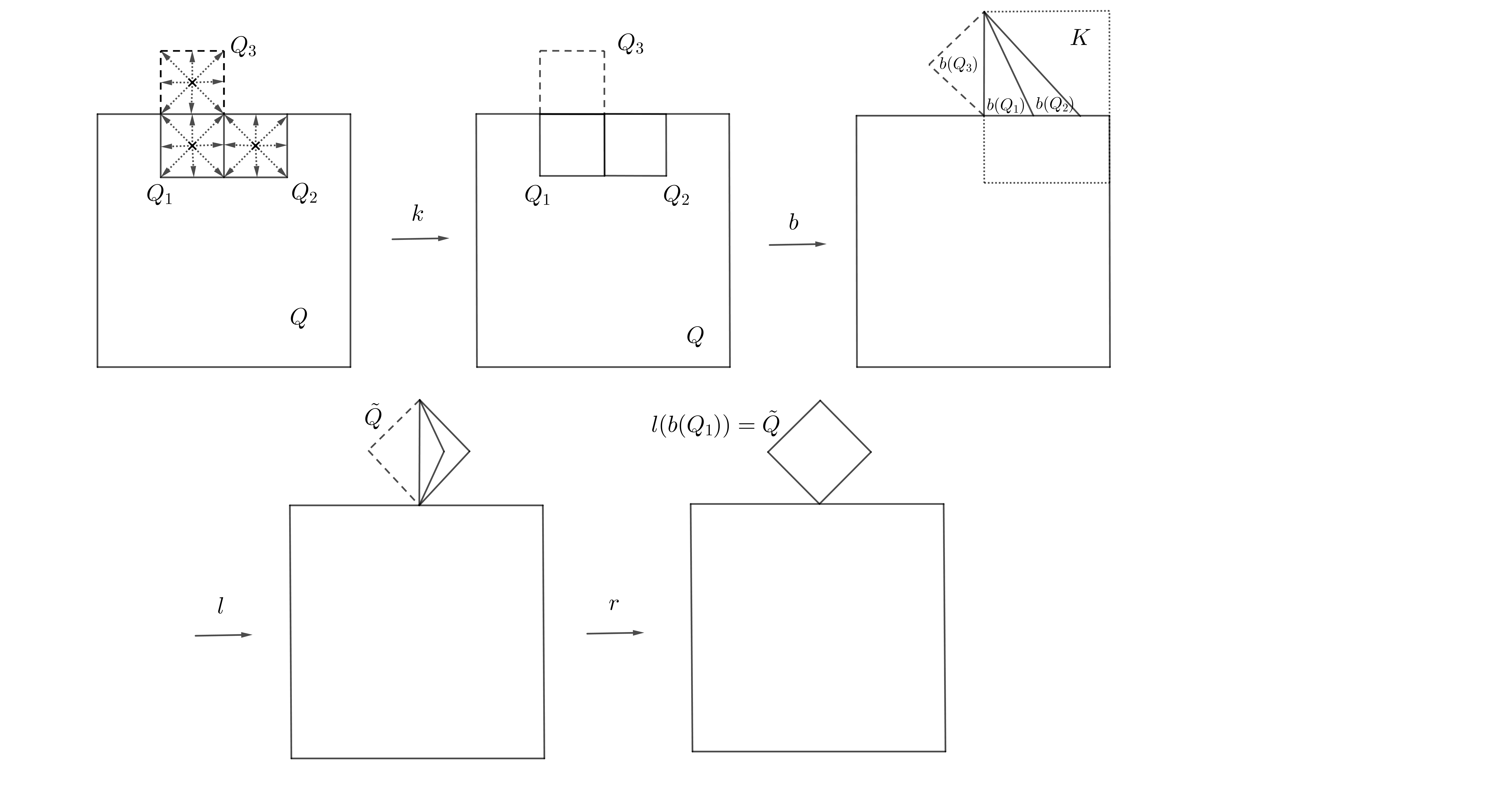}}
		\end{picture}
		\vskip -10pt	
		\caption{Construction of bad mapping in $(INV)$ which is not a mapping of finite distortion.}\label{noMFD}
	}
\end{figure}

The first mapping $k$ is the usual cavitation mapping in three cubes $Q_1$, $Q_2$, $Q_3$  and it sends $Q_i$ to $\partial Q_i$, i.e., the cavitation on the cube $Q(c,r)$ is given by formula
$$
k(x)=r\frac{x-c}{\|x-c\|_{\infty}}\text{ for }x\in Q(c,r),
$$
where $\|.\|_{\infty}$ denotes the maximum norm. It is easy to see that this mapping is continuous outside of $c$ and in a standard way we can show that it satisfies $k\in W^{1,p}([-2,2]^2,[-2,2]^2)$ for any $p<2$ (see e.g. the computation in \cite[Example 2.2]{HK}). As the mapping $r\circ l \circ b$ is Lipschitz, it is not difficult to see that our mapping $f\in W^{1,p}([-2,2]^2,[-2,2]^2)$ for any $p<2$. 

The construction of the second map $b$ can be seen from Fig. \ref{noMFD}. It maps $Q_1$, $Q_2$ and $Q_3$ onto the corresponding triangles, i.e., it puts $b(Q_1)$ and $b(Q_2)$ outside of $Q$. It is not difficult to see that such a mapping exists and we can construct it as a bi-Lipschitz map. It would be possible to give exact formulas but we believe that it is clear that such a mapping exists. 

The construction of the next mapping $l$ can be again seen from Fig. \ref{noMFD}. We squeeze triangles $b(Q_1)$ and $b(Q_2)$ in the $x$-direction so that their lower tip becomes sharp. The mapping is identity outside the box $K$ (see dotted rectangle in Fig. \ref{noMFD}) and it is of the form 
$$
l(x,y)=[a(x,y),y],
$$
where $a$ squeezes the horizontal segment in $b(Q_1)$ linearly onto the corresponding horizontal segment in $l\circ b(Q_1)$ and analogously it is linear on the part of the segment in $b(Q_2)$ and linear on the part of the segment $b(K\setminus (Q_1\cup Q_2))$. Again, it is easy to see how to construct this mapping and we do not present exact formulas.

The last mapping is identity outside of the (rotated) square $\tilde{Q}=l\circ b(Q_1\cup Q_2\cup Q_3)$ (see Fig. \ref{noMFD}). It again is of the form 
$$
r(x,y)=[\tilde{a}(x,y),y],
$$
it maps both $l\circ b(Q_2)$ and $l\circ b(Q_3)$ onto $\tilde{Q}$ so that it stretches each horizontal segment linearly. Now $l\circ b(Q_1)$ is also mapped onto $\tilde{Q}$ (so that it stretches each horizontal segment linearly) but now it reverses the orientation, i.e., it maps the left side of $\partial (l\circ b(Q_1))$ onto the right side of $\partial \tilde{Q}$ and vice versa, and its Jacobian is negative in $l\circ b(Q_1)$. This does  not matter it the end because our first mapping $k$ maps the whole $Q_1$ onto $\partial Q_1$ so $f(Q_1)=\partial \tilde{Q}$ and the Jacobian of $f$ is zero there.

It is easy to see that 
$$
f(Q_1)=\partial\tilde{Q},\ im_T (f,Q_1)=\tilde{Q}
\text{ and that } \deg(f,Q_1,y)=-1\text{ for }y\in \tilde{Q}
$$ 
as our mapping is changing the orientation there. Moreover, 
$$
f(Q_2)=\partial\tilde{Q},\ im_T (f,Q_2)=\tilde{Q}
\text{ and }\deg(f,Q_2,y)=1\text{ for }y\in \tilde{Q}. 
$$
Hence it is not difficult to see that 
$$
\deg(f,Q,y)=\deg(f,Q_1,y)+\deg(f,Q_2,y)=0\text{ for }y\in\tilde{Q}
$$
and $im_T (f,Q_1)=\tilde{Q}$ is not a subset of $im_T (f,Q)$. 

It remains to verify that this mapping satisfies the $(INV)$ condition. The full proof would require a lengthy case-by-case study. We give the details only for squares where the orientation is reversed and the reader can easily verify other cases. The orientation is reversed only for the square $Q_1$ so we should check other squares that are mapped to $im_T (f,Q_1)=\tilde{Q}$. Squares $Q_2$ and $Q_3$ are outside of $Q_1$ and  their images satisfy 
$$
f(Q_2)=f(Q_3)=\partial\tilde{Q}\subseteq \er^2\setminus im_T (f,Q_1)=\er^2\setminus\operatorname{int}(\tilde{Q}), 
$$
so the $(INV)$ condition holds there. On the other hand, $Q_1\subseteq Q$ and 
$$
f(Q_1)=\partial\tilde{Q}\subseteq im_T (f,Q)\cup f(\partial Q),
$$
as $\partial\tilde{Q}\subseteq f(\partial{Q})$, so the $(INV)$ condition holds there as well. Other cases can be discussed analogously.  
\end{proof}

The next example shows that without the condition that we have a mapping of finite distortion we cannot get that $h$ is Sobolev in Theorem \ref{main}, as we have only $BV$ regularity of the inverse there.

\prt{Lemma}
\begin{proclaim}\label{noMFDlemma}
	There is a Lipschitz mapping $f:[-1,4]\times[-1,2]\to[-1,3]\times[-1,2]$ which satisfies the $(INV)$ condition and equals to a homeomorphism close to the boundary, however, 
	$h$ is not Sobolev but only $BV$.
\end{proclaim}

\begin{proof}
For $[x,y]\in [0,3]\times[0,1]$ we define
$$
f(x,y)=
\begin{cases}
[x,y]&\text{ for }[x,y]\in [0,1]\times[0,1],\\
[1,y]&\text{ for }[x,y]\in [1,2]\times[0,1],\\
[x-1,y]&\text{ for }[x,y]\in [2,3]\times[0,1].\\
\end{cases}
$$
It is easy to check that this is a Lipschitz mapping which satisfies the $(INV)$ condition. Moreover, it is not difficult to find an extension to the whole $[-1,4]\times[-1,2]$ which is still Lipschitz and satisfies the $(INV)$ condition and which equals to a homeomorphism close to the boundary. 

Now it is easy to see that the natural inverse mapping $h$ (formally defined in the next section) satisfies
$$
h(x,y)=
\begin{cases}
[x,y]&\text{ for }[x,y]\in [0,1]\times[0,1],\\
[x+1,y]&\text{ for }[x,y]\in [1,2]\times[0,1]\\
\end{cases}
$$
and hence it has a jump on each segment $[0,2]\times\{t\}$ for every $t\in[0,1]$. It follows that $h\in BV\setminus W^{1,1}$. 

\end{proof}

\section{Definition of $h$ a.e. in $f_0(\Omega)$}\label{sec:h_def_ae}

It is not easy to define $h$ as we might have a set of zero measure $S\subseteq\Omega$, $|S|=0$, which is mapped onto a set of positive measure $T:=f(S)$, $|T|>0$, like in the Ponomarev example (see e.g. \cite[Theorem 4.10]{HK}). This mapping is a homeomorphism so it would be still possible to define $h$ as $f^{-1}$. However, we can compose it with a mapping which first squeezes a continuum $K_s$ to each point $s\in S$ and then maps it by the Ponomarev mapping to $t:=f(s)\in T$ (see \cite[Theorem 1.1]{BHM} for similar construction). Then naturally $f^{-1}_T(t)=K_s$ and it is not clear which point of $K_s$ we should choose when defining $h$ on a set of positive measure $T$. 

\subsection{Definition of $h$ in cavities} The point $x\in \Omega$ is called a cavity point if $|f_T(x)|>0$. By Theorem \ref{thm_embedded_ball} we know that for two different cavities $x_1, x_2\in \Omega$ we have
$$
|f_T(x_1)\cap f_T(x_2)|=0
$$
and hence it is easy to see that we can have at most countably many cavities. For cavity at $x$ we can use the $(INV+)$ condition for $B(x,1/n)$ and the complement of $B(x,1/k)$, $k<n$, and it is not difficult to deduce that for a.e. point in $y\in f_T(x)$ we know that there is no other point $a\in\Omega$ with $y\in f_T(a)$. Therefore for every cavity point $x$ and a.e. $y\in f_T(x)$ we can define 
$h(y)=x$ and this is well-defined. It is easy to see that for those $y\in f_T(x)$ that are moreover density points of $f_T(x)$ we have
\eqn{aaa}
$$
h(y)=\lim_{r\to 0+}\oint_{B(y,r)} h
$$
no matter how we define $h$ at other points (as $\Omega$ is bounded). We define the set of all images of cavities as
$$
C_{av}=\bigcup_{\{x\in\Omega:|f_T(x)|>0\}} f_T(x). 
$$

\subsection{Definition of $h$ at other points}

We would like to define $h$ at other points of $f_0(\Omega)$ by some analogy of \eqref{aaa}. We know from Lemma \ref{thm_injectivity} and Lemma \ref{diff2} that our mapping is $1-1$ a.e. on the set where $f$ is differentiable with $J_f>0$ and we denote
$$
G:=\{x\in \Omega:\ f\text{ is differentiable at }x, J_f(x)>0\}\setminus N,
$$
where $N$ is a set of measure zero such that $f$ is $1-1$ on $G$. Thus for every $y\in f(G)$ we define $h(y)=x$ where $x\in G$ satisfies $f(x)=y$. 

The following lemma tells us that around a.e. point in $f_0(\Omega)\setminus C_{av}$ we have plenty of points in $f(G)$ where $h$ is properly defined and hence it makes sense to define $h$ as
\eqn{limsup}
$$
h(y):=\limsup_{r\to 0+}\oint_{B(y,r)\cap f(G)}h\text{ for }y\in f_0(\Omega)\setminus C_{av},
$$ 
where the limes superior is taken component-wise.
In Lemma \ref{task5} below we show that this definition is correct for a.e. $y\in f(G)$, i.e., that it equals to $x\in G$ with $f(x)=y$ and that $\limsup$ could be replaced by $\lim$ in this case. 
The function $h$ is now defined a.e. in $f_0(\Omega)$, it serves as the correct analogy of the inverse function for us and in the next section we show that $h$ is Sobolev.

\prt{Lemma}
\begin{proclaim}\label{task3}
Let $\Omega\subseteq\rn$ be a domain and let $f_0:\Omega\to\rn$ be a homeomorphism. 
Let $f:\Omega\to\rn$ be a mapping of finite distortion in $W^{1,p}(\Omega, \mathbb{R}^n)$, $p>n-1$ (or $p\geq 1$ for $n=2$), that satisfies $(INV)$ and such that $f=f_0$ on a neighborhood of $\partial \Omega$. Then for a.e. $y\in f_0(\Omega)\setminus C_{av}$ and every $r>0$ we have
\eqn{goal}
$$
|f(G)\cap B(y,r)|>0.
$$
\end{proclaim}
\begin{proof}
By the change of variables formula \eqref{sub} and the fact that the Luzin $(N)$ condition holds for $f$ on the differentiability set of $f$ (see e.g. \cite[Corollary A.36 c)]{HK}) we obtain that 
$$
\bigl|f(\{x\in\Omega: J_f(x)=0, f\text{ is differentiable at }x\})\bigr|=0.
$$
Let us pick 
$$
y\in f_0(\Omega)\setminus (C_{av}\cup f(\{x\in\Omega: J_f(x)=0, f\text{ is differentiable at }x\})) 
$$
and $r>0$. We would like to prove \eqref{goal} for this $y$ and $r$. 

We first need to define some point $x_0\in\Omega$ whose every neighborhood is mapped close to $y$. We know that $f=f_0$ close to $\partial\Omega$ and that $f_0$ is a homeomorphism. Therefore
\eqn{one}
$$
\deg(f,\Omega,y)=\deg(f_0,\Omega,y)=1. 
$$
Without loss of generality we can assume that $\Omega$ is a union of finitely many cuboids (and that degree of $f$ is defined on the boundary of each such a cuboid as $f\in W^{1,p}\cap C$ on the boundary of these cuboids), as $f=f_0$ on some neighborhood of $\partial\Omega$ and by removing some part of this neighborhood we have a nicer domain $\Omega$ and still $f$ which equals to $f_0$ on some neighborhood of "new" $\partial\Omega$. Now we can divide $\Omega$ into many small cuboids (so that degree of $f$ is defined on the boundary of each such a cuboid) and by additivity of the degree and \eqref{one} we obtain that there is cuboid $R_0$ with 
$$
\text{ either }\deg(f,R_0,y)\neq 0\text{ or }y\in f(\partial R_0). 
$$
 We divide this $R_0$ into $2^n$ smaller cuboids and by additivity of the degree we either obtain that $\deg(f,R_1,y)\neq 0$ or that $y\in f(\partial R_1)$ for at least one of them. We continue by induction and we obtain a sequence of nested cuboids with sidelengths going to $0$, which gives us
$$
x_0:=\bigcap_{i=0}^{\infty} R_i\text{ so that }y\in im_T (f,R_i)\cup f(\partial R_1)\text{ for each }i.
$$

We claim that for every $i\in\en$ we have 
\eqn{positive}
$$
|R_i\cap G|>0. 
$$
Assume for contrary that this is not the case. As our mapping is differentiable a.e. we obtain that 
$J_f(x)=0$ a.e. in $R_i$. Since $f$ is a mapping of finite distortion we obtain that $|Df(x)|=0$ a.e. in $R_i$. Using the ACL condition in all $n$ directions we now obtain that a.e. point of $R_i$ is mapped to a single point $z$. As  $y\in im_T (f,R_i)\cap f(\partial R_1)$ we obtain that $z=y$ and thus 
$$
y\in f\bigl(\{x\in\Omega: J_f(x)=0, f\text{ is differentiable at }x\}\bigr)
$$ 
which is a contradiction to our choice of $y$. 

\begin{figure}[h t p]
	\vskip 200pt
	{\begin{picture}(0.0,0.0) 
			\put(-150.2,0.2){\includegraphics[width=0.90\textwidth]{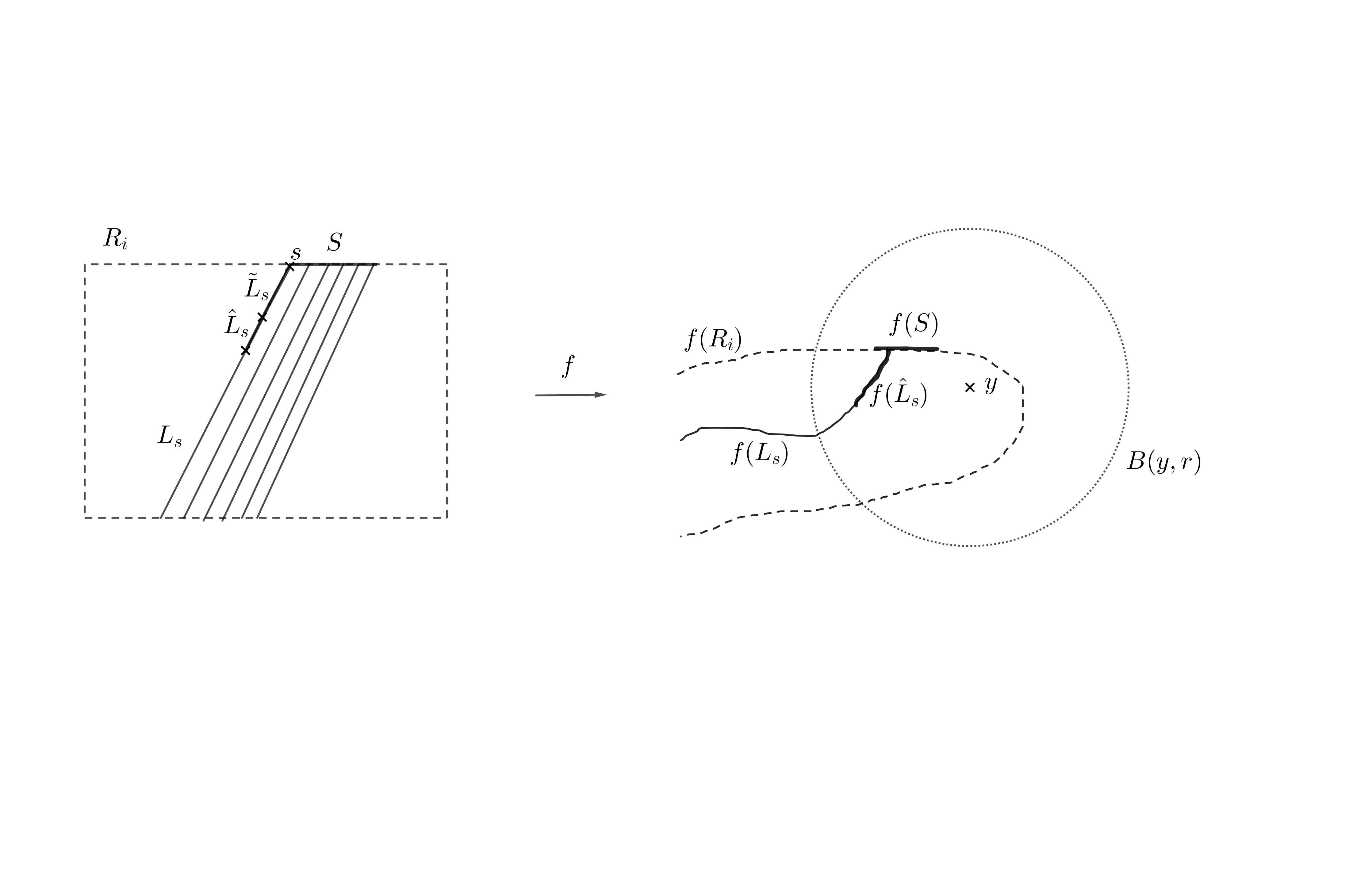}}
		\end{picture}
		\vskip -90pt	
		\caption{Finding of big part of $f(G)$ in $B(y,r)$ for the validity of \eqref{goal} for $n=2$.}\label{ahoj}
	}
\end{figure}

Now let us return to proving \eqref{goal}. If $y\in f(\partial R_i)$ then we can clearly find $(n-1)$-dimensional disc  
\eqn{aha} 
$$
S\subseteq \partial R_i\text{ such that }
f(S)\subseteq B(y,r)  
$$
(recall that $f$ is actually continuous on $\partial R_i$, so $f(S)$ is well-defined). 
We claim that also in the case $\deg(f,R_i,y)\neq 0$ for all $i$ we can find $i_0$ big enough so that for every $i\geq i_0$ we find $(n-1)$-dimensional disc $S$ such that \eqref{aha} holds, see Fig. \ref{ahoj}. Indeed, otherwise we would have $B(y,r)\subseteq f_T(x_0)$, but $y$ was chosen not to be a point of cavity as $y\notin C_{av}$. 

Without loss of generality $S$ lies only on one side of $R_i$. By using the ACL condition we can find 
$S_0\subseteq S$ with $\haus^{n-1}(S_0)>0$ and many parallel segments $L_s$ inside $R_i$ so that $f$ is absolutely continuous on $L_s$, $s\in S_0$, and each segment $L_s$ ends in $s\in S_0$, see Fig. \ref{ahoj}. We know \eqref{positive} and thus we can assume without loss of generality that 
$$
|G\cap \bigcup_{s\in S_0}L_s|>0
$$
as we can choose a direction of the parallel segments and there is definitely some direction in which 
$S$ "can see" a part of $R_i\cap G$ which is big enough. As our mapping is differentiable a.e. and of finite distortion, we can assume that $f$ is differentiable $\haus^1$-a.e. and it satisfies $J_f(x)=0\Rightarrow |Df(x)|=0$ $\haus^1$-a.e. on $L_s$ for each $s\in S_0$. 
Now let us pick one $s\in S_0$ with 
\eqn{aaaa}
$$
\haus^1(L_s\cap G)>0.
$$ 
We find a maximal segment $\tilde{L}_s\subseteq L_s$ which end in $s\in S_0$ such that $J_f=0$ a.e. on 
$\tilde{L}_s$, see Fig. \ref{ahoj}. As $f$ is a mapping of finite distortion we obtain that $|Df(x)|=0$ a.e. on this segment, and as $f$ is absolutely continuous there, we obtain that $f(\tilde{L}_s)=s$. 
We can now find a segment $\hat{L}_s\subseteq L_s\setminus \tilde{L}_s$ next to 
$\tilde{L}_s$ (they have one common endpoint) so that by continuity of $f$ and 
$$
f(s)\in f(S)\subseteq B(y,r)
$$
we have
$$
f(\hat{L}_s)\subseteq B(y,r).
$$
 As we know that $f$ is differentiable a.e. on $L_s$, $J_f$ is not $0$ a.e. on $\hat{L}_s$ (by the definition of $\tilde{L}_s$) and that \eqref{aaaa} holds, we obtain that $\haus^1(\hat{L}_s\cap G)>0$ and hence the set
 $$
 A:=\bigcup_{s\in S_0}\hat{L}_s\cap G\text{ satisfies }|A|>0\text{ and }f(A)\subseteq B(y,r). 
 $$
 Since $J_f\neq 0$ on $A$ by the definition of $G$, we obtain by the change of variables formula that $|f(A)|>0$ and \eqref{goal} follows. 
\end{proof}

\subsection{The definition of $h$ is correct at $f(G)$} 
The following lemma tells us that the definition of $h$ is correct, i.e., at a.e. point of $f(G)$
the natural definition of $h$ coincides with \eqref{limsup}.

\prt{Lemma}
\begin{proclaim}\label{task5}
Let $\Omega\subseteq\rn$ be a domain and let $f_0:\Omega\to\rn$ be a homeomorphism. 
Let $f:\Omega\to\rn$ be a mapping of finite distortion in $W^{1,p}(\Omega, \rn)$, $p>n-1$ (or $p\geq 1$ for $n=2$), that satisfies $(INV)$ and such that $f=f_0$ on a neighborhood of $\partial \Omega$. For $x\in G$ we set $y=f(x)$ and we have
\eqn{eee}
$$
x=\lim_{r\to 0+}\oint_{B(y,r)\cap f(G)}h. 
$$ 
\end{proclaim}
\begin{proof}
Let us pick $x\in G$. We know that $J_f(x)>0$ and by differentiability
\eqn{diff}
$$
\lim_{a\to x}\frac{f(a)-f(x)-Df(x)(a-x)}{|a-x|}=0.
$$
We denote $y=f(x)$ and for $a\in G$ we set $z=f(a)$ (which is well-defined and $h(z)=a$ as $a\in G$). Recall that by Lemma \ref{task3} we know that we have plenty of points $z$ like that around $y$.  

Note that \eqref{diff} especially implies (as $J_f(x)>0$) that
$$
|z-y|=|f(a)-f(x)|\approx|Df(x)(a-x)|\approx |a-x|=|h(z)-h(y)|. 
$$
From \eqref{diff} we should thus get (as $J_f(x)>0$) that 
$$
\begin{aligned}
0&=\lim_{z\to y,\ z\in f(G)}\frac{(Df(x))^{-1}(z-y)-(h(z)-h(y))}{|h(z)-h(y)|}\\
&=\lim_{z\to y,\ z\in f(G)}\frac{h(z)-h(y)-(Df(x))^{-1}(z-y)}{|z-y|}.\\
\end{aligned}
$$
Our conclusion \eqref{eee} now follows easily.
\end{proof}

\section{Regularity of the inverse}\label{sec:reg_of_inv}

In this section we will prove that the partial inverse map defined in previous section satisfies $h\in W_{\loc}^{1,1} (\Omega', \rn)$ and it is of finite distortion. Our starting point is to prove the following auxiliary inequality.
\subsection{A key estimate}

\begin{thm}\label{thm:key_estimate}Let $\Omega, \Omega'\subseteq \rn$ be domains and let $f_0:\Omega\to \Omega'$ be a homeomorphism. 
	Let $p>n-1$ (or $p\ge n-1$ when $n=2$).  Suppose that $f\in W^{1,p}(\Omega,\rn)$ is a mapping of finite distortion which satisfies the $(INV)$ condition and $f=f_0$ on a neighborhood of $\partial \Omega$. Then
	\begin{equation}\label{eq:key_estimate}
		\diam (f^{-1} (B_r)) \le C(p,n) r^{1-n} \abs{f^{-1} (B_{2r})}^\frac{p-n+1}{p} \left( \int_{f^{-1} (B_{2r})} \abs{Df}^p\right)^\frac{n-1}{p}
	\end{equation}
	for all balls $B_r = B(y,r) \subseteq \rn$ with $B_{3r}\subseteq\Omega'$.
\end{thm}
Recall that $f^{-1} (B_r)$ denotes the preimage of the ball $B_r$. By Lemma \ref{lemma_strong} and Theorem \ref{thm_inv+} we can assume that $f$ satisfies strong $(INV+)$.  

We may and do assume that  $ \diam f^{-1}(B_r)=d$ and the set $\overline{f^{-1}(B_r)}$ contains the origin and the point $(d, 0, \dots , 0)$. For $0\le t \le d$, we denote 
\[L_t = \{x=(x_1, \dots , x_n) \in \Omega  \colon x_1=t\} \, . \]

\begin{lemma}\label{lem:connectness}
	For a.e. $t\in (0,d)$ we have $L_t \cap f^{-1}(B_r) \not= \emptyset$. 
\end{lemma}
\begin{proof}
	Let's proceed by contradiction and assume that there exists a set   $U \subseteq (0,d)$ such that  $\mathcal{H}^1(U) > 0$ and $L_t \cap f^{-1}(B_r) = \emptyset$ for every $t\in U$.  Since a.e. point of $U$ is a density point of $U$ and $\mathcal{H}^1(U) > 0$, we can pick $u\in U$ of density 1. Now we can find $u_k, v_k\in U$ with $u_k \nearrow u$ and $v_k\searrow u$ and moreover we can assume that the hyperplanes 
	$\{x_1=u_k\}\cap \Omega$ and $\{x_1=v_k\}\cap \Omega$ are so nice that we can use the $(INV+)$ condition for sets whose boundary contains part of these hyperplanes.  
	
	Now we find sets $U_k, V_k$ with Lipschitz boundary (see Fig. \ref{ahaaa}) for which we can use the $(INV+)$ condition, $U_k\subseteq U_{k+1}$, $V_k\subseteq V_{k+1}$, 
	\eqn{all}
	$$
	\bigcup_{k=1}^{\infty} U_k=\Omega\cap\{x_1<u\},\quad \bigcup_{k=1}^{\infty} V_k=\Omega\cap\{x_1>u\}\text{ and }(\partial U_k\cup\partial V_k)\cap f^{-1}(B_r)=\emptyset. 
	$$
	
		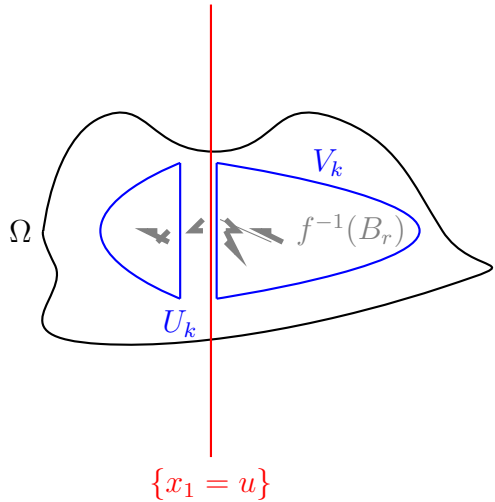
\begin{figure}[h t p]
		\begin{tikzpicture}[scale=1.50]
			\draw[thick, black] plot[smooth, tension=1] coordinates {(3.5,0)(4,1)(5,0.7)(6.05,1)(7,0)(7,-0.5)(4,-1)(3.6,-0.3)(3.5,0)};
			\node at (3.3,0) {$\Omega$};

			\node[gray] at (6.2,0) {$f^{-1}(B_r)$};
			\draw[thick, blue] (5.02, 0.6) -- (5.02, -0.6) ;
			\draw[thick, blue] plot[smooth, tension=1] coordinates {(5.02, 0.6)(6.8,0) (5.02, -0.6)};
			\node[blue] at (6,0.6) {$V_k$};
			
			\draw[thick, red] (4.97, 2) -- (4.97, -2) node[below] {$\{x_1=u\}$};
			\draw[thick, blue] (4.7, 0.6) -- (4.7, -0.6) node[below] {$U_k$};
			\draw[thick, blue] plot[smooth, tension=1] coordinates {(4.7, 0.6)(4,0) (4.7, -0.6)};

			\draw[gray] plot[smooth, tension=1] coordinates {(5.5,-0.1)(5.05,0.1)};
			\draw[line width=0.7mm, gray] plot coordinates {(4.5,-0.1)(4.6,0)};
			\draw[line width=0.7mm, gray] plot coordinates {(5.1,-0.1)(5.2,-0.2)(5.1,0)};
			\draw[line width=0.7mm, gray] plot coordinates {(5.6,-0.1)(5.4,0)(5.5,0)};
			\draw[line width=0.7mm, gray] plot coordinates {(4.6,-0.1)(4.4,0)(4.5,0)};
			\draw[line width=0.7mm, gray] plot coordinates {(4.9,0.1)(4.8,0)(4.9,0)};
			\draw[line width=0.7mm, gray] plot coordinates {(5.1,0.1)(5.2,0)(5.1,0)};
			
		\end{tikzpicture}
			\caption{The construction of $U_k$ and $V_k$.}\label{ahaaa}
		\end{figure}	
	
	To construct $U_k$ we just follow the hyperplane $\{x_1=u_k\}$ close to $\partial \Omega$ where $f=f_0$ (and thus it is $1-1$ and it does not intersect $f^{-1}(B_r)$) and we construct a Lipschitz curve close to $\partial \Omega$ so we get a Lipschitz domain $U_k$ such that $ f^{-1}(B_r)\cap\{x_1\leq u_k\}\subseteq$ $U_k \subseteq \{x_1\leq u_k\}$ and for which we can use $(INV+)$. To construct $U_{k+1}$ we proceed similarly, we just keep $U_k$ nested in $U_{k+1}$. To construct $V_k$ we similarly use hyperplanes $\{x_1=v_k\}$ and we construct $ f^{-1}(B_r)\cap\{x_1\geq v_k\}\subseteq$ $V_k\subseteq \{x_1\geq v_k\}$. 
	
By Theorem \ref{thm_embedded_ball} (which holds also for more general shapes than balls with similar proof) we obtain that $\{U_k\}$ have nested topological images and that topological images of $U_k$ and $V_k$ are disjoint. Since $(\partial U_k\cup\partial V_k)\cap f^{-1}(B_r)=\emptyset$, we obtain from the strong $(INV+)$ condition that (from a certain index $k_0$ up)
	$$
	\emptyset\neq f(U_k)\cap B_r\subseteq im_T(f,U_k)\text{ and }\emptyset\neq f(V_k)\cap B_r\subseteq im_T(f,V_k).
	$$
	By \eqref{all} and $u\in U$ we obtain that
	$$
	B_r=\Bigl(\bigcup_k  im_T(f, U_k)\cap B_r\Bigr) \cup \Bigl(\bigcup_k  im_T(f,V_k)\cap B_r\Bigr),
	$$
	but the connected set $B_r$ cannot be written as a union of two disjoint nonempty open sets, which gives us a contradiction.

\end{proof}

\begin{remark}\label{rem:connectness}
The same result also holds if we replace $f^{-1}(B_r)$ by $f^{-1}(y)$, where $y\in\Omega'$ -- in the final step of the proof we get that 
$$
y\in im_T(f, U_k)\cap im_T(f,V_k) = \emptyset,
$$
obtaining the contradiction. Actually the result holds for any connected compact set.
\end{remark}

\begin{proof} [Proof of Theorem~\ref{thm:key_estimate}]
	We may assume that $d=\diam f^{-1}(B_r)>0$ because otherwise there is nothing to prove. By the Fubini theorem $f|_{L_t} \in W^{1,p} (L_t , \rn)$ for almost every $t\in (0,d)$.  Next we fix $t\in (0,d)$ so that $L_t \cap f^{-1}(B_r) \not= \emptyset$ and  $f|_{L_t} \in W^{1,p} (L_t , \rn)$.   Since $f = f_0$ near the boundary of $\Omega$ and $B_{3r}\subseteq\Omega'$, $f(L_t) \not\subseteq B_{2r}$. Also we know that $f(L_t) \cap \partial B_r \neq \emptyset$. Mapping $f$ is continuous on $L_t$, so $f(L_t)$ is a nice connected (locally) compact set which has to span from the boundary of $B_{2r}$ into $B_r$ as on Fig. \ref{obrazek_lajna}.

	\begin{figure}[h t p]
		\begin{tikzpicture}
			\draw (0,0) circle (1);
			\draw (0,0) circle (2);
			
			\draw[thick, red] plot[smooth, tension=1] coordinates {(2.5,0) (1,1) (0.5,0.5) (0,1) (-1,0)(-1,2)(-1.5,1)(-2,2)};
	  	\node[red] at (2.5,0.6) {$f(L_t)$};
			
			\node at (-1.6,-0.5) {$B_{2r}$};
			\node at (-0.6,-0.5) {$B_r$};    
			\draw[<-, thick] (-2.5, 0) -- (-4, 0) node[midway, above] {$f$};
			
			\node at (-7,0) {$f^{-1}(B_r)$};
			\draw[thick, red] (-5, 2) -- (-5, -2) node[below] {$L_t$};
			
			\draw[gray] plot[smooth, tension=1] coordinates {(-5.5,-0.1)(-5.05,0.1)};
			\draw[line width=0.7mm, gray] plot coordinates {(-4.5,-0.1)(-4.6,0)};
			\draw[line width=0.7mm, gray] plot coordinates {(-5.1,-0.1)(-5.2,-0.2)(-5.1,0)};
			\draw[line width=0.7mm, gray] plot coordinates {(-5.6,-0.1)(-5.4,0)(-5.5,0)};
			\draw[line width=0.7mm, gray] plot coordinates {(-4.6,-0.1)(-4.4,0)(-4.5,0)};
			\draw[line width=0.7mm, gray] plot coordinates {(-4.9,0.1)(-4.8,0)(-4.9,0)};
			\draw[line width=0.7mm, gray] plot coordinates {(-5.1,0.1)(-5.2,0)(-5.1,0)};
			
		\end{tikzpicture}
			\caption{The set $f(L_t)$ goes from $S_{2r}$ inside $B_r$ and back.}\label{obrazek_lajna}
\end{figure}
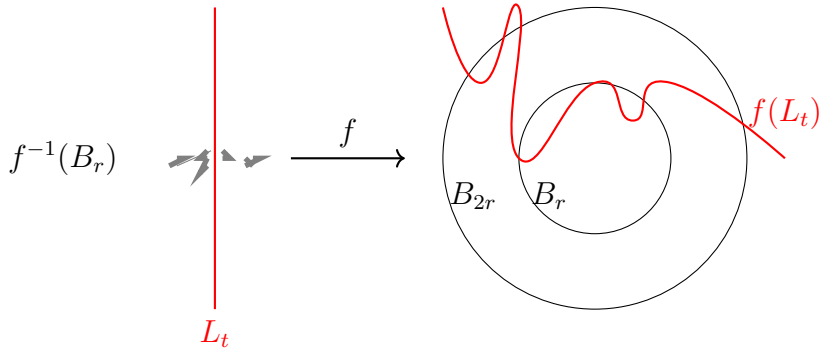	

	Fix a point $x(t) \in L_t \cap f^{-1} (B_r)$ and  choose a radius $s(t)>0$ such that a $(n-1)$-dimensional ball 
	$$
	B^{n-1} (x(t), s(t)) \subseteq L_t\cap f^{-1} (B_{2r})\text{ and }\partial B^{n-1} (x(t), s(t)) \cap \partial  f^{-1} (B_{2r})\neq \emptyset.  
	$$  
	To simplify the notation we denote $B(t) = B^{n-1} (x(t), s(t))$. Now, the mapping $f$ is continuous on $L_t$ and Lemma~\ref{lem:soboembed} gives 
		\[r \le C(p,n) s(t)^{1- \frac{n-1}{p}} \left( \int_{B(t)} \abs{Df}^p\right)^\frac{1}{p}.\]
	This estimate can be equivalently written as
	\[s(t)^{n-p-1} \le [C(p,n)]^p r^{-p} \int_{B(t)} \abs{Df}^p.\]
	By Lemma~\ref{lem:connectness} this inequality holds for a.e. $t\in (0,d)$ and, therefore, we can integrate it from $0$ to $d$ to obtain
	\begin{equation}\label{eq:almostthere}
		\int_0^d s(t)^{n-p-1}\; dt \le [C(p,n)]^p r^{-p} \int_{f^{-1}(B_{2r})} \abs{Df}^p.
	\end{equation}
	Here we also used the fact that $B(t) \subseteq L_t\cap f^{-1} (B_{2r})$. 
	
	If $n=2$ and $p=1$, then the claimed estimate~\eqref{eq:key_estimate} follows immediately from~\eqref{eq:almostthere}.
	
	If $p>n-1$, then for estimating the integral on the left-hand
	side, we set a positive number $\alpha = (n-1)(p-n+1)/p$ and apply the H\"older inequality to obtain
	\[
	\begin{split}
		d^\frac{p}{n-1} &= \left(\int_0^d s(t)^{-\alpha} s(t)^\alpha \; dt\right)^\frac{p}{n-1} \le \left(\int_0^d s^{-\alpha \frac{p}{n-1}} \right) \left(\int_0^d s^{\alpha \frac{p}{p-n+1}} \right)^\frac{p-n+1}{n-1} \\
		&=  \left(\int_0^d s^{n-p-1} \right) \left(\int_0^d s^{n-1} \right)^\frac{p-n+1}{n-1} \\
		& \le \left(\int_0^d s^{n-p-1} \right) \abs{f^{-1} (B_{2r})}^\frac{p-n+1}{n-1}.
	\end{split}
	\]
	Combining this with~\eqref{eq:almostthere} the claimed estimate~\eqref{eq:key_estimate}  also follows  when $p>n-1$. 
\end{proof}
Recall that
$$
G=\{x\in \Omega \colon  J_f(x)>0 \text{ and }f \text{ is differentiable and 1-1}\}. 
$$
Under the  assumptions of Theorem~\ref{thm:key_estimate}, the estimate~\eqref{eq:key_estimate} gets also the following
form:
\begin{equation}\label{eq:key_estimate_ver2}
	\diam (f^{-1} (B_r)) \le C(p,n) r^{1-n} \abs{f^{-1} (B_{2r})}^\frac{p-n+1}{p} \left( \int_{B_{2r}} \Phi(y)\right)^\frac{n-1}{p}
\end{equation}
where
\eqn{defPhi}
$$
\Phi (y) :=\frac{\abs{(Df)(h(y))}^p}{J_f(h(y))} \chi_{f(G)}(y). 
$$

Indeed,  first since our mapping $f$ has finite distortion we have
\begin{align*}
	\int_{f^{-1} (B_{2r})} \abs{Df(x)}^p \, d x &= \int_{\{ x\in f^{-1} (B_{2r})\colon J_f(x)>0 \}} \abs{Df(x)}^p \, d x + \int_{\{ x\in f^{-1} (B_{2r})\colon J_f(x)=0 \}} \abs{Df(x)}^p \,d x \\
	&= \int_{\{ x\in f^{-1} (B_{2r})\colon J_f(x)>0 \}} \abs{Df(x)}^p \, d x \, . 
\end{align*}
Second, we know from Lemma \ref{thm_injectivity} and Lemma \ref{diff2} that our mapping is $1-1$ a.e. on the set where $J_f>0$ and $f$ is differentiable a.e. in $\Omega$. Thus
\[ 
\int_{f^{-1} (B_{2r})} \abs{Df(x)}^p \, d x =  \int_{f^{-1} (B_{2r})\cap G} \abs{Df(x)}^p \,d x  =  \int_{f^{-1} (B_{2r})\cap G} \frac{\abs{Df(x)}^p}{J_f(x)} J_f(x) \, d x . 
\]
Recall that for every $y\in f(G)$ we defined $h(y)=x$, where $x\in G$ satisfies $f(x)=y$.  Now, since our mapping $f$ is differentiable in $x\in G$, we can use the change of variables formula (see e.g. \cite[Corollary A.36 (c)]{HK}), and we have
\[\int_{f^{-1} (B_{2r})} \frac{\abs{Df(x)}^p}{J_f(x)} J_f(x) \chi_G (x) \, d x = \int_{B_{2r}} \Phi (y) \, d y \, .  \]
Note that this identity also shows that $\Phi \in L^1 ({B_{2r}})$.

We claim that the estimate~\eqref{eq:key_estimate_ver2} holds also for $h$, i.e., 
\begin{equation}\label{eq:key_estimate_ver3_aux}
	\abs{h(y_1)-h(y_2)} \le C(p,n) r^{1-n} \abs{f^{-1} (B_{2r})}^\frac{p-n+1}{p} \left( \int_{B_{2r}} \Phi(y)\right)^\frac{n-1}{p}
\end{equation}
for almost every $y_1, y_2 \in B_r$. Note that the mapping $h\in L^\infty(f_0 (\Omega))$ and so far $h$ is only defined almost everywhere in $f_0 (\Omega)$, see Section~\ref{sec:h_def_ae}. After proving the estimate~\eqref{eq:key_estimate_ver3_aux} we will show which representative of $h$ satisfies this estimate for every $y_1, y_2 \in B_r$.

There are three cases that we distinguish based on the definition of $h$ at those points. 
If $y\in f(G)\cap B_r$ then $h(y)\in f^{-1}(y)$ and the estimate obviously holds on $h(B_r\cap f(G))$ as it is the same as \eqref{eq:key_estimate_ver2}. If $y\in B_r\cap(f_0(\Omega)\setminus C_{av})$ then we can use Lemma \ref{task3} and \eqref{limsup} to obtain that 
$$
h(y)=\limsup_{r\to 0+}\oint_{B(y,r)\cap f(G)}h.
$$
As the estimate~\eqref{eq:key_estimate_ver2} holds for $f^{-1}$ it is not difficult to see that it holds also for integral averages (of values for which the estimate holds). It remains to consider the case $y\in B_r\cap C_{av}$. If $B_r\subseteq f_T(x)$ then there is nothing to prove as $h=x$ a.e. on $B_r$, so we may assume that 
$$
\haus^{n-1}\left(B_r\cap \partial f_T(x)\right)>0. 
$$
We know that 
$$
f_T(x)=\bigcap_{\delta>0, \delta \notin N_{h(y)}} E(f^*,B(h(y),\delta))
$$ 
and that $f_T(x)\subseteq E(f^*,B(h(y),\delta))$ for a.e. $\delta$ (where the topological image is well-defined). 
We choose such a $\delta$ (such that $f\in W^{1,p}(\partial B(h(y),\delta),\rn)$ so it is continuous there) and we obtain that each closed connected set $f(\partial B(h(y),\delta))$ has to encircle $f_T(x)$. Set $S_\delta:=\partial B(h(y),\delta)$.
As $\haus^{n-1}\left(B_r\cap \partial f_T(x)\right)>0$ we may assume that $\delta$ is small enough so that 
\eqn{rrr}
$$
\haus^{n-1}\bigl(f(S_\delta\cap B_r)\bigr)>0.
$$
Without loss of generality we also assume that $\haus^{n-1}$-a.e. point in $S_\delta$ is a point of differentiability of $f$, that $f|S_\delta$ is of finite distortion and that $\haus^{n-1}$-a.e. point in $S\cap\{J_f>0\}$ belongs to $G$. 
We claim that we can find a point 
\eqn{rrrr}
$$
\tilde x\in G\cap S\cap f^{-1}(B_r).
$$  
If this is not the case, then for $\haus^{n-1}$-a.e. point in $S_\delta\cap f^{-1}(B_r)$ we have $J_f=0$ and hence $Df=0$, which implies that also all $(n-1)\times (n-1)$ minors of $Df$ are zero. 
Since $f\in W^{1,p}(S_\delta,\rn)$ with $p>n-1$, we obtain that images of $\haus^{n-1}$-zero subsets of $S_\delta$ have zero $\haus^{n-1}$-measure (see e.g. \cite{MM}). The same conclusion holds for $n=2$ by the ACL condition on the circle. It follows that the $(n-1)$-dimensional change of variables formula (analogy of Theorem \ref{subst} in those curved coordinates)  holds for $f$ on $S_\delta$ and since 
$Df=0$ a.e. on $S_\delta\cap f^{-1}(B_r)$ we obtain that 
$$
\haus^{n-1}\bigl(f(S_\delta\cap B_r)\bigr)=0.
$$
That contradicts \eqref{rrr} and we get that there exists $\tilde x$ as in \eqref{rrrr}. 
Since $\tilde x\in G$ we obtain that  \eqref{eq:key_estimate_ver3_aux} holds for this $\tilde x=h(\tilde y) \in f^{-1}(\tilde y)$ for some $\tilde y$ and $|\tilde x-h(y)|<\delta$. We can choose $\delta$ as small as we wish, so the estimate holds also for our $y$.

Next we will show that the average representative of $h$,
$$
\tilde{h}(y):= \limsup_{s\to 0+}\oint_{B(y,s)}h(z)\, d z,
$$
satisfies~\eqref{eq:key_estimate_ver3_aux} for every $y_1, y_2 \in B_r$, i.e.,
\begin{equation}\label{eq:key_estimate_ver3}
	\diam \tilde{h} (B_r) \le C(p,n) r^{1-n} \abs{f^{-1} (B_{2r})}^\frac{p-n+1}{p} \left( \int_{B_{2r}} \Phi(y)\right)^\frac{n-1}{p}.
\end{equation}
Indeed, let $y\in B_r$. Choose $y_1 \in B_r$ so that~\eqref{eq:key_estimate_ver3_aux} holds at $y_1$.  Then
$$
 \abs{\oint_{B(y,s)}h(z) \, d z -h(y_1)} \le C(p,n) r^{1-n} \abs{f^{-1} (B_{2r})}^\frac{p-n+1}{p} \left( \int_{B_{2r}} \Phi(y)\right)^\frac{n-1}{p}
 $$
for $0<s<\dist(y, \partial B_r)$. Letting $s \to 0+$ we have
$$
\abs{\tilde{h}(y) -h(y_1)} \le C(p,n) r^{1-n} \abs{f^{-1} (B_{2r})}^\frac{p-n+1}{p} \left( \int_{B_{2r}} \Phi(y)\right)^\frac{n-1}{p}.
$$
Since $y \in B_r$ is arbitrary, the claimed estimate~\eqref{eq:key_estimate_ver3} follows.

\begin{corollary}\label{dusl}
Let $p\in (n-1,n]$ (or $[n-1,n]$ if $n=2$). Then a.e. point in $f(\Omega)\setminus f(G)$ has exactly one preimage. Moreover, for a.e. $y\in\Omega'$ we have $y\in f_T(h(y))$. 
\end{corollary}
\begin{proof}
Let $y\in f(\Omega)\setminus f(G)$ be a Lebesgue point of $\Phi\in L^1$. From \eqref{eq:key_estimate_ver2} we get (for $r<1$)
\begin{align*}
	\diam (f^{-1} (y))  \leq \diam (f^{-1} (B(y,r))) &\leq \frac{\diam (f^{-1} (B(y,r)))}{r^{(n/p-1)(n-1)}}\\
	&\leq C(p,n, \Omega') \left(r^{-n}\int_{B(y,2r)} \Phi(y)\right)^{\frac{n-1}{p}},
\end{align*}
and as our $\Phi$ contains $\chi_{f(G)}$ (see \eqref{defPhi}) we obtain that the right-hand side tends to $0$ as $r\to 0$ since $y$ is a Lebesgue point of $\Phi$ and $y\notin f(G)$.

Let $y\in\Omega'$. We want to show that $y\in f_T(h(y))$. If $y\in C_{av}$, $h(y)$ is defined as $x$ such that $y\in f_T(x)$, and if $y\in f(G)$, then $h(y)\in f^{-1}(y)\cap G$. The remaining case follows from \eqref{eq:key_estimate_ver3}: for $r>0$ we find $y'\in B(y,r)\cap G$. Then we have

\begin{align*}
	\diam (f^{-1} (y)\cup h(y)) &\leq \diam (f^{-1} (y)\cup (f|G)^{-1} (y')) + \diam ((f|G)^{-1} (y')\cup h(y))\\
	& \leq \diam (f^{-1} (B(y,r))) + \diam (h(B(y,r))) \\
	& \leq  C(p,n, \Omega') \left(r^{-n}\int_{B(y,2r)} \Phi(y)\right)^{\frac{n-1}{p}}
\end{align*}
and so for $y\in \Omega'\setminus(f(G)\cup C_{av})$ which is a Lebesgue point of $\Phi$ we have $f^{-1}(y)=\{h(y)\}$.
\end{proof}

\begin{proof}[Proof of Theorem \ref{main} and Theorem \ref{main2}]\phantom{a}\newline 
  {\it Step 1:} {\bf  $\tilde{h}$ is differentiable almost everywhere} 
  
  By the Lebesgue differentiation theorem~\cite[Theorem 23.5]{Vab}, the set function $E \mapsto \abs{f^{-1} (E)}$ has finite derivative $\mu' (y)$ at almost every $y$; that is,
 $$
 \mu' (y) = \lim_{r \to 0} \frac{\abs{f^{-1} \big(B(y,r)\big)}}{\abs{B(y,r)}} < \infty \quad \textnormal{ a.e. } y \in \Omega' \, . 
 $$
 By the Rademacher-Stepanov theorem, $\tilde{h}$ is differentiable
 almost everywhere in $\Omega' $ provided that
 \begin{equation}\label{eq:rade_step}
 	\limsup_{r\to 0} \frac{\diam \tilde{h}(B_r(y)) }{r} < \infty  \quad \textnormal{ a.e. } y \in \Omega' \, . 
 \end{equation}
 Now fix a point $y\in \Omega'$ such that 
 $$
 \mu' (y) < \infty  \quad \textnormal{and} \quad \limsup_{r\to 0} \oint_{B(y,r)} \Phi (z) \, d z = \Phi (y) < \infty \, .  
 $$
 Notice that almost every point in $\Omega'$ satisfies these conditions. By~\eqref{eq:key_estimate_ver3} we have
 \begin{equation}\label{eq:diff_est}
 	\limsup_{r\to 0} \frac{\diam \tilde{h}(B(y,r))}{r}  \le C(p,n) [\mu'(y)]^\frac{p-n+1}{p} [\Phi (y)]^\frac{n-1}{p} < \infty \, .
 \end{equation}
 Hence $\tilde{h}$ is differentiable almost everywhere in $\Omega'$. For this pointwise derivative we claim that
 \begin{equation}\label{eq:Dh_int}
 	\abs{Dh(y)} \in L^1_{\loc} (\Omega') \, . 
 \end{equation}
 Indeed, this follows from~\eqref{eq:diff_est} after applying the H\"older inequality because the
 Lebesgue theorem~\cite[Theorem 23.5]{Vab} also gives that
 $$
  \int_E \mu' (y) d y \le  \abs{f^{-1} (E)} \, .  
  $$

  {\it Step 2:} { \bf  $|Dh|=0$ a.e. outside of $f(G)$}
  
  We know that $h$ is differentiable a.e., so let's take $y$ a point of differentiability of $h$. By \eqref{eq:diff_est} we have
$$
|Dh(y)| \leq \limsup_{r \to 0+} \frac{\mathrm{diam}(h(B(y, r)))}{r} \leq C_p(n)[\mu'(y)]^{\frac{p-n+1}{p}}[\Phi(y)]^\frac{n-1}{p}<\infty
$$
for almost every $y$.
Since
$$
\Phi(y) = \frac{|Df(h(y))|^p}{J_f(h(y))}\chi_{f(G)}(y),
$$
the right-hand side is 0 outside of $f(G)$. Therefore $|Dh(y)| = 0$ a.e. outside of $f(G)$.

 {\it Step 3:} { \bf  Sobolev regularity $h\in W^{1,1} (\Omega' , \rn)$} 
 
 For that we need only to
 show that $\tilde{h}$ is ACL, i.e., absolutely continuous on almost all lines parallel to the
 coordinate axes, since the estimate~\eqref{eq:Dh_int} guarantees the local integrability of the
 partial derivatives. For the ACL property, fix an open cube $Q$ such 
 that $\overline{Q} \subseteq \Omega^\prime.$ It suffices to show 
 that $\tilde{h}$ is ACL in this cube. By symmetry, it suffices to consider 
 line segments parallel to the $x_n$-axis. Assume that 
 $Q=  Q_0 \times  I$. Next, for each Borel set 
 $ E\subseteq  Q_0$ we define
 \begin{equation}
 	\eta ( E) = \abs{f^{-1} ( E \times I)}
 \end{equation}
 By the Lebesgue differentiation theorem, $\eta$ 
 has finite 
 derivative $\eta^\prime (y)$ for $\haus^{n-1}$-a.e. $y \in Q_0$, that is,
 $$
 \eta^\prime (y) = \limsup_{r\to 0} \frac{\eta (B^{n-1} (y,r))}{r^{n-1}} < \infty \, . 
 $$
 Denote by $V_0$ the $\haus^{n-1}$-zero set in $Q_0$ where $\eta$ does not exist or is infinity.
 Also, since  $\Phi \in L^1_{\loc} (\Omega')$, there exists a $\haus^{n-1}$-zero set in $Q_0$ (denoted by $V_1$) such that
 $$
 \int_I \Phi (y,t) \, d t < \infty \qquad \textnormal{for all } y\notin V_1 \, . 
 $$
 
 Let $\mathcal A (I) = \mathcal A $ be set of finite unions of closed intervals in $I$ whose end-points are
 rational numbers. Clearly, this set is countable. For all $J\in \mathcal A$ we define a
 function $\phi_J \colon Q_0 \to \mathbb R$ such that
 $$
 \phi_J (y) = \int_J  \Phi (y,t) \, d t  \qquad \textnormal{when } y\in Q_0 \setminus V_1 \, . 
 $$
 Then $\phi_J \in L^1 (Q_0)$. Thus $\haus^{n-1}$-a.e. $y\in Q_0$ is a Lebesgue point for
 $\phi_J$. Denote by $V_J$ the $\haus^{n-1}$-zero set of non-Lebesgue points. Now
 \[V= V_0 \cup V_1 \cup \Bigl( \bigcup_{J\in \mathcal A} V_J \Bigr) \]
 has $\haus^{n-1}$-measure zero.
 
 Fix $y\in Q_0 \setminus V$. We will prove that $h$ is absolutely
 continuous on the segment $\{y\}\times I$, which proves the
 claim. Let $\{\triangle_i\}_{i=1}^\ell$, $\triangle_i=[a_i,b_i]$, be a union of $\ell$
 closed intervals on ${I}$ whose interiors are mutually disjoint  and whose endpoints
 are rational numbers.
 
 Fix a natural number $\beta$ such that 
 $$
 \frac{1}{\beta}<\frac{
 	\min\big\{\mbox{dist}\big(\cup {\Delta}_i, \partial {I}\big),
 	b_1-a_1,\dots,b_\ell-a_\ell \big\}}{100}.
 $$
  Now a standard covering 
 argument \cite[Lemma 31.1]{Vab} gives us a number $\delta$ depending on 
 $\beta$ such
 that for $0<r<\delta$ we have a covering
 $B_1^i,\dots,B_{\kappa_i}^i$ of $\triangle_i$ which has
 the properties
 \begin{itemize}
 	\item $\mbox{diam}(B_j^i)=2r$, 
 	\item overlapping of $B_j^i$ is at most $4$, 
 	\item
 	$B_j^i \subseteq \frac{1}{\beta}$-neighborhood of $\triangle_i$.
 \end{itemize}
 Define $n$-dimensional balls ${\mathbb B}^i_j={\mathbb B}(x_j^i,r)$ such that ${\mathbb B}^i_j\cap \{y\}\times I=\{y\}\times B_j^i.$ Without loss of generality we can assume that ${\mathbb B}(x_j^i,4r)\subseteq \Omega'$. 
 To shorten our notation we 
 write ${\mathbb A}^i_j= f^{-1}({\mathbb B}(x_j^i,2r))$. Now, applying the estimate~\eqref{eq:key_estimate_ver3} and the H\"older inequality we have
 \[
 \begin{split}
 	& \sum_{i=1}^\ell \left|\tilde{h}(y,a_i)- \tilde{h}(y,b_i)\right|  \le
 	\sum_{i=1}^\ell \sum_{j=1}^{\kappa_i} \mbox{diam} (\tilde{h} ( {\mathbb B}_j^i)) \\
 	&\le  C(p,n) \sum_{i=1}^\ell \sum_{j=1}^{\kappa_i}  r^{1-n} \abs{\mathbb A_j^i}^\frac{p-n+1}{p} \left( \int_{4\mathbb B_j^i} \Phi(y) \, d y\right)^\frac{n-1}{p}
 	\\
 	&= C(p,n) \sum_{i=1}^\ell \sum_{j=1}^{\kappa_i}  \left( \frac{\abs{\mathbb A_j^i}}{r^{n-1}} \right)^\frac{p-n+1}{p} \left( \frac{1}{r^{n-1}} \int_{4\mathbb B_j^i} \Phi(y) \, d y\right)^\frac{n-1}{p}  \\
 	&\le C(p,n) \left( \frac{1}{r^{n-1}}\sum_{i=1}^\ell \sum_{j=1}^{\kappa_i}   \abs{\mathbb A_j^i}  \right)^\frac{p-n+1}{p}  \left(\frac{1}{r^{n-1}}  \sum_{i=1}^\ell \sum_{j=1}^{\kappa_i} \int_{4\mathbb B_j^i} \Phi(y) \, d y \right)^\frac{n-1}{p}. 
 \end{split}
 \]
 From the geometry of our covering it follows that the overlapping of 
 $$
 \sum_{i=1}^\ell \sum_{j=1}^{\kappa_i}   \chi_{4 \mathbb B_j^i} \le C  .
 $$
 Therefore, we have
 \begin{equation}
 	\sum_{i=1}^\ell \sum_{j=1}^{\kappa_i}   \abs{\mathbb A_j^i} \le C \eta \big(B^{n-1} (y,4r) \big)
 \end{equation}
 and
 \begin{equation}
 	\sum_{i=1}^\ell \sum_{j=1}^{\kappa_i} \int_{4\mathbb B_j^i} \Phi(y) \, d y \le C
 	\int_{B_{4r} \times J^\beta } \Phi(y) \, d y\,  ,
 \end{equation}
 where $J^\beta = \bigcup_{i=1}^\ell [a_i-1/\beta, b_i+1/\beta] $.
 By taking the upper limit when $r \rightarrow 0$ and then letting $\beta \rightarrow \infty,$ we conclude that
 \begin{equation}\label{eq:sobo_int}
 	\sum_{i=1}^\ell \left| \tilde{h} (y,a_i) - \tilde{h} (y,b_i)\right|  \le
 	C_p(n) [\eta^\prime (y)]^\frac{p-n+1}{p} \left(\int_{\cup_{i=1}^\ell [b_i-a_i] }  \Phi(y, t) \, d t\right)^\frac{n-1}{p}
 \end{equation}
This establishes that $\tilde{h}$ is absolutely continuous on the line segment $y \times I$, as desired. Consequently, we have completed the proof that $h \in W_{\mathrm{loc}}^{1,1} (\Omega', \mathbb{R}^n)$. Observe that it suffices to consider rational endpoints in~\eqref{eq:sobo_int}, since one can first show (using analogous estimates) that $\tilde{h}$ is continuous on $y \times I$ for almost every \(y \in Q_0\). 

Henceforth, we will identify \(h\) with \(\tilde{h}\), unless explicitly stated otherwise. 

 {\it Step 4:} {\bf The mapping $h$ is of finite distortion} 
 
By Step 2 we know that $Dh(y)=0$ for a.e. $y\notin f(G)$ so there is nothing to prove. 
 Second, we assume that $y\in f(G)$. Now since $h$ is differentiable in a.e. $y\in \Omega'$, we know from the proof of Lemma~\ref{task5} that
 \begin{equation}\label{eq_inv_diff_matrix}
 Dh (y) = [Df \big(h(y) \big)]^{-1} \quad \textnormal{ for almost every } y\in f(G) \, .
 \end{equation}
In particular, for the corresponding $x\in G$ with $f(x)=y$ we have
 \begin{equation}\label{eq_jacb_inv}
 J_h(y) = \frac{1}{J_f (h(y))}=\frac{1}{J_f(x)}  \qquad \textnormal{ for a.e. } y\in f(G) \, . 
 \end{equation}  
Since $0<J_f(x)<\infty$ for $x\in G$ we know that $h$ is of finite distortion as desired. 

Theorem \ref{main} and Theorem \ref{main2} now follow once we show the $(INV)$ condition for $h$ for $n=2$ in the next section. 
\end{proof}

\section{(INV) for $h$ in the case $n=2$}\label{sec:2d_h_inv}

In the planar case, the regularity $h\in W^{1,1}$ is enough to talk about the $(INV)$ condition for $h$. (This is not true in higher dimensions, as there we would need at least $h\in W^{1,n-1}$.)

\prt{Lemma}
\begin{proclaim}\label{lines_intersection}
Let $M$ be a segment and let $g=(g_1,g_2)\in W^{1,1}(M,\er^2)$. Then for $\haus^1$-a.e. $y$ the set $g(M)$ intersects the line 
$$
L_y:=\{(x,y), x\in \mathbb{R}\}
$$
in only finitely many points and with $D g_2\neq 0$ at those points. 
\end{proclaim}

\begin{proof}
Since $g$ is absolutely continuous on $M$, it maps null sets to null sets and hence 
$$
\haus^1\bigl(g(\{s: D g_2(s)\text{ does not exist}\})\bigr)=0.
$$
Moreover, by the change of variables formula
$$
\haus^1\left(g_2(\{s:D g_2(s)= 0\}\right)=0.
$$
It follows that we avoid both these set with $L_y$ for $\haus^1$-a.e. $y$. 
The finiteness of $g(M)\cap L_y$ follows from the Eilenberg inequality (see e.g. \cite{EH}) where we choose $f(x,y)=y$ and $s=t=1$. It gives
$$
\int_\mathbb{R} \haus^0(g(M)\cap (\mathbb{R}\times \{y\}))\leq C \haus^1(g(M))<\infty,
$$
so the integrand on the left-hand side must be finite for a.e. $y$.
\end{proof}

\prt{Theorem}
\begin{proclaim}\label{thm_h_inv}
Let $\Omega\subseteq\mathbb{R}^2$ be a domain and let $f_0:\Omega\to\mathbb{R}^2$ be a homeomorphism. 
Let $f:\Omega\to\mathbb{R}^2$ be a mapping of finite distortion in $W^{1,p}$, $p\geq 1$, that satisfies the $(INV)$ condition
and such that $f=f_0$ on a neighbourhood of $\partial \Omega$. Then $h$ satisfies $(INV)$.
\end{proclaim}

\begin{proof}
First, it is enough to consider only balls $B$ such that $h \in W^{1, 1}(\partial B, \mathbb{R}^2)$, 
\eqn{rr}
$$
J_h = 0\implies |Dh|=0 \text{ } \haus^{1}\text{-a.e. on } \partial B
$$
 and (due to the line above and Step 2 in the proof of Theorem \ref{main} and \ref{main2}) 
\eqn{rr2}
$$
J_h(y) > 0\implies y\in f(G) \text{ } \haus^{1}\text{-a.e. on } \partial B.
$$

Without loss of generality we can assume that our ball $B_r=B(0,r)$.
Let's assume for contradiction that there exists $R$ a set of radii, $\haus^1(R)>0$, and a set $A_r \subseteq B_r$ such that 
$$
|A_r| > 0 \text{ and } h(A_r) \cap E(h, B_r) = \emptyset \text{ for every } r\in R.
$$

 {\it Step 1:} {\bf  We show that we can assume $A_r \subseteq f(G)$}

We can find a direction such that segments connecting $A_r$ and $S_r = \partial B_r$, on which $h$ is $W^{1,1}$ and of finite distortion, form a set of positive measure. For such segment $L_a$, we know that $h(L_a)$ connects $h(S_r)$ with $h(a)$, $a\in A_r$, which lies outside of $E(h, B_r)$.

Since $E(h, B_r)$ and $h(a)$ are disjoint compact sets, 
$$
\haus^1 (h(L_a) \setminus E(h, B_r)) > 0.
$$
 Denote $L'_a$ the component of $L_a \setminus h^{-1}(E(h, B))$ which contains $a\in A_r$.
Since $h$ is absolutely continuous on $L'_a$, it cannot happen that $|Dh| = 0$ $\haus^1$-a.e. on $L'_a$. Therefore there exists $\tilde{L}_a \subseteq L'_a$ such that $\haus^1(\tilde{L}_a) > 0$ and $J_h > 0$ on $L'_a$ $\haus^1$-a.e. Without loss of generality $\tilde{L}_a \subseteq f(G)$ upto a set of $\haus^1$-measure zero. Since 
$$
h(\tilde{L}_a) \cap E(h, B_r) = \emptyset \text{ and }L_a\subseteq B_r,
$$
 we can pick "new" $A_r$ as 
 $$
 \cup_{a\in A_r}\tilde{L}_a\cap f(G).
 $$

 {\it Step 2:} {\bf We can have a dyadic cube $Q$ which is mostly contained in $A_r$ for a lot of $r\in R$}
 
We can assume that 
there exists $\delta>0$ such that for all $r\in R$ we have $|A_r|>\delta$.  We can fix $r_1\in R$ a point of density of $R$ and find $r_0<r_1$ such that $|B_{r_1}\setminus B_{r_0}|<\delta/2$. 

By using Lebesgue density theorem for each $r\in R\cap [r_0, r_1]$ we can find $k_r$ such that for every $k\geq k_r$ we have a dyadic cube $Q\subseteq B_{r_0}$ of sidelength $2^{-k}$ with 
$$
|Q\cap A_r|>|Q|/2.
$$
Now we have divided a set of positive measure $R\cap [r_0, r_1]$ into countably many sets according to $k_r$ and hence one of them has to have a positive $\haus^1$ measure.  We call this new set of radii $R_0$ and for each $r\in R_0$ we have a dyadic cube $Q_r$ of sidelength $2^{-k}$ (where $k$ is the same for all $r\in R_0$). As there are only finitely many dyadic cubes of sidelength $2^{-k}$ inside $B_{r_0}$ there has to be one of them which works for a subset of $R_0$ of positive $\haus^1$ measure. We call this set of good radii $R_1$ and we work only with it further on. Set $Q':= Q\cap f(G)$.

{\it Step 3:} {\bf Finding a good line $L$}

Since $f$ is a homeomorphism near $\partial f(\Omega)$, we can find $C \subseteq f(G)$:
$$
|C| > 0, \quad C \cap B_r = \emptyset, \quad h(C) \cap E(h, B_r) = \emptyset
$$
for all $r\in R_1$.

Since $h|_{f(G)} = (f|_G)^{-1}$ and $f$ satisfies the Luzin $(N)$ condition on $G$, we have that $|h(A_r)| > 0$ and $|h(C)| > 0$.

We denote the circle as $S_r:=\partial B_r$. By $h\in W^{1,1}(S_r,\er^2)$, \eqref{rr} and \eqref{rr2} we can assume that for any $r \in R_1$ we have $h(S_r) \subseteq G$ up to a set of $\haus^1$-measure zero.

So now we pick a segment $L$ connecting $h(a), a\in Q'$, and $h(c), c\in C$, which satisfies
\begin{itemize}
\item[(a)] $f$ is absolutely continuous on $L$ and $|Df| = 0$ $\haus^1$-a.e. on $L \setminus G$,
\item[(b)] $L$ and $h(S_r)$ intersect in only finitely many points, which all lie in $G$ and have only one preimage under $h|_{S_r}$, for $r \in R_L \subseteq R_1$, where $\haus^1(R_L) > 0$, 
\item[(c)] we can moreover assume that whenever $L$ intersects $h(S_r)$, every neighbourhood of such point on $L$ contains points both in $im_T(h, B_r)$ and outside of $E(h, B_r)$ (i.e., $L$ goes in or out of $im_T(h, B_r)$, it doesn't just touch the boundary),
\item[(d)] the endpoint $h(a)$ lies in $h(A_r)$ for $r\in R_3\subseteq R_L$ where $\haus^1(R_3) > 0$.
\end{itemize}

The property (a) is immediate. We get (b) from the fact that $\haus^1(h(S_r \setminus f(G))) = 0$, so we can avoid the set with a lot of lines, and from Lemma \ref{lines_intersection} (when using the polar coordinates). To prove (c) we assume without loss of generality that $L$ is parallel to the $x$ axis and we use again Lemma \ref{lines_intersection}. In the intersection point $w\in L\cap h(S_r)$ we know that $Dg_2(w)\neq 0$ and hence we are not only touching $g(S_r)$ there, but we have to go through (i.e. we are in the situation of green lines on Fig. \ref{fig5} and not of the red line). Moreover, 
$w$ has only one preimage under $h|_{S_r}$ due to (b), so it cannot happen that two parts of $h(S_r)$ touch in $w$ (as in the pinched point on the red line on Fig. \ref{fig6}). 
Since (a), (b) and (c) hold for a.e. such segment, using the Fubini theorem we obtain (d).

\begin{center}
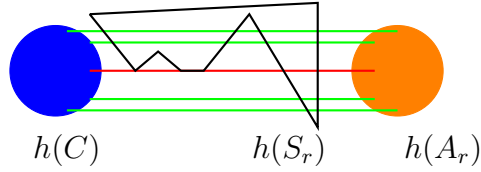
\begin{figure}[h t p]
		\begin{tikzpicture}[scale=1.50]

	\draw[orange, fill=orange] (1,0) circle (0.4);
	\draw [blue, fill=blue] (-2,0) circle (0.4);	
	
				\node at (0.05,-0.7) {$h(S_r)$};   
				\node at (1.4,-0.7) {$h(A_r)$};   
			\node at (-1.9,-0.7) {$h(C)$};    

\draw[thick, red] (-1.7,0)--(0.8,0);
\draw[thick, green] (-1.7,0.25)--(0.8,0.25);
\draw[thick, green] (-1.9,0.35)--(1,0.35);
\draw[thick, green] (-1.7,-0.25)--(0.8,-0.25);
\draw[thick, green] (-1.9,-0.35)--(1,-0.35);
	
	\draw[thick] (-1.7,0.5) -- (-1.3,0) -- (-1.1, 0.17) -- (-0.9,0)--(-0.7,0)--(-0.3,0.5) -- (0.3,-0.5) --(0.3,0.6)--(-1.7,0.5);

		\end{tikzpicture}
		\caption{The green horizontal lines intersect $h(S_r)$ in finitely many points with a positive derivative in the direction perpendicular to those lines. The red line (in the middle) contains both points where the derivative does not exist and where it is zero.}\label{fig5}
		\end{figure}
\end{center}

As neither $h(a)$ nor $h(c)$ are in $E(h, B_r)$ and $L$ does not intersect $h(S_r)$ in a "pinched" point (i.e., where the $h|_{S_r}$-preimage of such point is not connected), the set 
$$
D_r:=L \cap h(S_r)= L \cap f^{-1}(S_r) \cap G
$$ 
must have an even number of elements, see Fig. \ref{fig5} and Fig. \ref{fig6} where green segments intersect $h(S_r)$ in even number of points. 

\begin{center}
\begin{figure}[h t p]
		\begin{tikzpicture}[scale=1.50]
			\draw[fill=gray] (1,0) circle (2);
			\draw[fill=white] (0,0) circle (1);
						
			\node at (1.7,-0.5) {$E(h,B)$};

	\draw[orange, fill=orange] (0,0) circle (0.4);
	\draw [blue, fill=blue] (-2,0) circle (0.4);	
	
				\node at (0.05,-0.7) {$h(A)$};    
			\node at (-1.9,-0.7) {$h(C)$};    

\draw[thick, red] (-1.7,0)--(-0.2,0);
\draw[thick, green] (-1.7,0.25)--(-0.2,0.25);
\draw[thick, green] (-1.9,0.35)--(0,0.35);
\draw[thick, green] (-1.7,-0.25)--(-0.2,-0.25);
\draw[thick, green] (-1.9,-0.35)--(0,-0.35);
	
		\end{tikzpicture}
		\caption{The green horizontal lines intersect $h(S_r)$ in an even number of points. The red line (in the middle) intersects it only in one, pinched point.}\label{fig6}
		\end{figure}
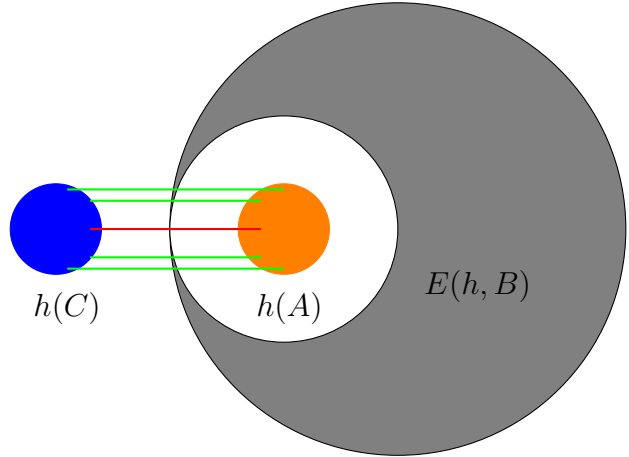
\end{center}

{\it Step 4:} {\bf $f$-image of $L$}

Now we look at the $f$-images: $f(h(a)) = a$, $f(h(c)) = c$, so $f(L)$ goes from inside of $B_r$ outside. Again, for a lot of radii ($r\in R_3$ with $\haus^1(R_3)>0$), we know from Lemma \ref{lines_intersection} (used with polar coordinates) that $f(L)$ intersects $S_r$ always with non-zero derivative in the radial direction, and only finitely many times. Since we start in $B_r$ and end outside, we need to cross $S_r$ odd-many times. Note that in general $f^{-1}(S_r)$ might be bigger than $h(S_r)$ as some points in $S_r$ might have more preimages under $f$. 
We denote 
$$
N_r := L \cap f^{-1}(S_r) \setminus G.
$$
Now the set $D_r\cup N_r$ is the set of all points on $L$ whose $f$-images intersect $S_r$. We know $\haus^0(D_r)$ is even and $\haus^0(D_r \cup N_r)$ is odd, therefore $N_r \neq \emptyset$. That means that there are some points on $L$ which do belong to $f^{-1}(S_r)$ but not to $h(S_r)$. Denote
$$
N=\bigcup_{r\in R_3} N_r.
$$

That gives us $\haus^1(\Pi_r(f(N))) > 0$, so $\haus^1(f(N)) > 0$. However, $N \subseteq L \setminus G$, so $|Df| = 0$ a.e. on $N$, and $f$ is absolutely continuous  on $L$, so $\haus^{1}(f(N)) = 0$, which gives us the desired contradiction.

{\it Step 5:} {\bf The case $A_r \subseteq f(\Omega) \setminus B_r$, $h(A_r) \subseteq im_T(h, B_r)$}

This case is done analogously.  The differences are that in Step 2 $Q\subseteq  f(\Omega)\setminus B_{r_1}$ and in Step 4 $\haus^0(D_r)$ is odd and $\haus^0(D_r \cup N_r)$ is even, as we are connecting a point $h(a)$ lying inside $im_T(h, B_r)$ and a point $h(c)$ lying outside of $E(h, B_r)$.

\end{proof}

Since the situation between $f$ and $h$ is now entirely symmetrical for $n=2$ we obtain the following analogue of Corollary \ref{dusl}. 

\begin{corollary}
Let $n=2$ and $f\in W^{1,1}_{\loc}(\Omega, \rn)$. Then for a.e. $x\in\Omega$ we have $x\in h_T(f(x))$. 
\end{corollary}

\section{Proof of Theorems~\ref{thm_cor_1} and~\ref{thm_cor_2}}\label{sec:conf_inner}
\begin{proof}
By Step 2 of the proof of Theorem \ref{main} we know that for a.e. $y\notin f(G)$ we have $|Dh(y)|=0$. Therefore,
$$
\int_{\Omega'} |Dh(y)|^n \, dy = \int_{\Omega'} |Dh(y)|^n  \chi_{f(G)} (y)\, d y .
$$

Theorems~\ref{main} and~\ref{main2} give that $h$ is differentiable almost everywhere in $\Omega'$ and is of finite distortion.  Therefore,
$$
\int_{\Omega'} |Dh(y)|^n \, dy = \int_{\mathcal A} |Dh(y)|^n  \chi_{f(G)} (y)\, dy,
$$
where 
$$
\mathcal A = \{ y \in \Omega'  \colon J_h(y) > 0 \text{ and } h \text{ is differentiable at } y \} \, . 
$$
In particular,
$$
\int_{\Omega'} |Dh(y)|^n \, dy = \int_{\mathcal A} \frac{|Dh(y)|^n}{J_h(y)}  \chi_{f(G)} (y) J_h(y)\, dy,
$$
and after changing variables, applying Theorem~\ref{subst} for $h$, we have
$$
\int_{\Omega'} |Dh(y)|^n \, dy \le \int_{h(\mathcal A )} \frac{|(Dh)(f(x))|^n}{J_h(f(x))}  \chi_{G} (x)  dx = \int_{h(\mathcal A )} |(Dh)(f(x))|^n J_f(x)  \chi_{G} (x) \, dx \, . 
$$
The last identity follows from~\eqref{eq_jacb_inv}. We have just proved that
\begin{equation}\label{eq_ine_1}
\int_{\Omega'} |Dh(y)|^n \, dy \le  \int_{\Omega} |(Dh)(f(x))|^n J_f(x)  \chi_{G} (x) \, dx \, . 
\end{equation}
On the other hand, the change of variable  formula (Theorem~\ref{subst}) applied to $f$ gives
\begin{equation}\label{eq_ine_2}
\int_{\Omega'} |Dh(y)|^n \, dy \ge  \int_{\Omega} |(Dh)(f(x))|^n J_f(x)  \chi_{G} (x) \, dx \, . 
\end{equation}
Combining \eqref{eq_ine_1} with~\eqref{eq_ine_2}, we find that
$$
\int_{\Omega'} |Dh(y)|^n \, dy =  \int_{\Omega} |(Dh)(f(x))|^n J_f(x)  \chi_{G} (x) \, dx . 
$$
Since $(Dh) (f(x))= [Df(x)]^{-1}$ for a.e. $x\in G$ by~\eqref{eq_inv_diff_matrix}, we have
$$
\int_{\Omega'} |Dh(y)|^n \, dy =  \int_{\Omega} |[Df(x)]^{-1}|^n J_f(x)  \chi_{G} (x) \, dx . 
$$
The Cramer rule gives
$$
\int_{\Omega} |[Df(x)]^{-1}|^n J_f(x)  \chi_{G} (x) \, dx  = \int_{\Omega '} \frac{\lvert D^\sharp f(x) \rvert^n}{J_f(x)^{n-1}} \chi_{G} (x) \, dx
= \int_{\Omega'} K_I(x, f) \chi_{G} (x)\, dx
$$
and the desired identity follows.

For $h\in W^{1,n}(\Omega',\Omega)$ we can mimic the proof of Lemma \ref{lines_intersection} to get an analogous result for an $(n-1)$-dimensional ball $M$, $g\in W^{1,n}(M,\rn)$ and for $\haus^{n-1}$-a.e. $y$ and $L_y=\{(x,y_1,\dots,y_{n-1}), x\in \mathbb{R}\}$.
Then, following the lines of the proof of Theorem \ref{thm_h_inv}, we obtain that $h$ satisfies the $(INV)$ condition.
\end{proof}


\begin{thebibliography}{GKL}



\bibitem{A}
\by{\name{Adams}{D. R.}}
\paper{A note on Choquet integrals with respect to Hausdorff capacity}
\jour{In: Function Spaces and Applications (Lund, 1986). Lecture Notes in Mathematics 1302, Springer}
\pages{115--124}
\yr{1988}
\endpaper



\bibitem{AIMb}
\by{\name{Astala}{K.}, \name{Iwaniec}{T.} and \name{Martin}{G.}}
\book{Elliptic partial differential equations and quasiconformal mappings in the plane}
\publ{Princeton University Press}
\yr{2009}
\endbook


\bibitem{AIMO}
\by{\name{Astala}{K.}, \name{Iwaniec}{T.}, \name{Martin}{G. J.} and \name{Onninen}{J.}}
\paper{Extremal mappings of finite distortion}
\jour{Proc. London Math. Soc. (3)}
\vol{91}
\yr{2005}
\pages{655--702}
\endpaper



\bibitem{Ba}
\by{\name{Ball}{J.}}
\paper{Global invertibility of Sobolev functions and the
interpenetration of matter}
\jour{Proc. Roy. Soc. Edinburgh Sect. A}
\vol{88\nom 3--4}
\pages{315--328}
\yr{1981}
\endpaper

\bibitem{BHMC}
\by{\name{Barchiesi}{M.}, \name{Henao}{D.}, \name{Mora-Corral}{C.} and \name{Rodiac}{R.}}
\paper{Harmonic dipoles and the relaxation of the neo-Hookean energy in 3D elasticity}
\jour{Arch. Ration. Mech. Anal.}
\vol{247} 
\yr{2023}
\pages{Paper No. 70, 46 pp}
\endpaper


\bibitem{BHMCR2}
\by{\name{Barchiesi}{M.}, \name{Henao}{D.}, \name{Mora-Corral}{C.} and \name{Rodiac}{R.}}
\paper{On the lack of compactness problem in the axisymmetric neo-Hookena model}
\jour{Forum Math. Sigma}
\vol{12} 
\yr{2024}
\pages{Paper No. e26, 70 pp}
\endpaper



\bibitem{BP}
\by{\name{Bauman}{P.} and \name{Phillips}{D.}}
\paper{Univalent minimizers of polyconvex functionals in 2 dimensions}
\jour{Arch. Rational Mech. Anal.}
\vol{126}
\yr{1994}
\pages{161--181}
\endpaper


\bibitem{BPO1}
\by{\name{Bauman}{P.}, \name{Phillips}{D.} and \name{Owen}{N.}}
\paper{Maximum principles and a priori estimates for an incompressible material in nonlinear elasticity}
\jour{Comm. Partial Differential Equations}
\vol{17}
\yr{1992}
\pages{1185--1212}
\endpaper


\bibitem{BHM}
\by{\name{Bouchala}{O.}, \name{Hencl}{S.} and \name{Molchanova}{A.}}
\paper{Injectivity almost everywhere for weak limits of Sobolev homeomorphisms}
\jour{J. Funct. Anal.}
\vol{279} 
\yr{2020}
\pages{article 108658}
\endpaper



\bibitem{CN}
\by{\name{Ciarlet}{P. G.}, \name{Ne\v{c}as}{J.}}
\paper{Injectivity and self-contact in nonlinear elasticity}
\jour{Arch. Ration. Mech. Anal.}
\vol{97}
\pages{171--188}
\yr{1987}
\endpaper

\bibitem{CDL}
\by{\name{Conti}{S.}, \name{De Lellis}{C.}}
\paper{Some remarks on the theory of elasticity for compressible
Neohookean materials}
\jour{ Ann. Sc. Norm. Super. Pisa Cl. Sci.}
\vol{2}
\pages{521--549}
\yr{2003}
\endpaper

\bibitem{CHM}
  \by{\name{Cs\"ornyei}{M.}, \name{Hencl}{S.} and \name{Mal\'y}{J.}}
  \paper{Homeomorphisms in the Sobolev space $W^{1,n-1}$}
  \jour{J. Reine Angew. Math}
  \vol{644}
  \pages{221--235}
  \yr{2010}
  \endpaper

\bibitem{DPP}
\by{\name{De Philippis}{G.} and \name{Pratelli}{A.}}
\paper{The closure of planar diffeomorphisms in Sobolev spaces}
\jour{Ann. Inst. H. Poincare Anal. Non Lineaire}
\vol{37}
\pages{181--224}
\yr{2020}
\endpaper

 \bibitem{DHM}
 \by{\name{Dole\v{z}alov\'a}{A.}, \name{Hencl}{S.}, \name{Mal\'y}{J.}}
 \paper{Weak limit of homeomorphisms in $W^{1,n-1}$ and $(INV)$ condition}
 \jour{Arch. Rational Mech. Anal.}
 \vol{247}
 \pages{Article No. 80, 54 pp}
 \yr{2023}
 \endpaper

 \bibitem{DHMo}
 \by{\name{Dole\v{z}alov\'a}{A.}, \name{Hencl}{S.}, \name{Molchanova}{A.}}
 \paper{Weak limit of homeomorphisms in $W^{1,n-1}$: invertibility and lower semicontinuity of energy}
 \jour{ESAIM: Control Optim. Calc. Var.}
 \vol{30}
 \pages{article 37}
 \yr{2024}
 \endpaper




\bibitem{EG}
\by{\name{L.C.}{Evans} and \name{R.F.}{Gariepy}}
\book{Measure theory and Fine Properties of Functions}
\publ{Studies in Advanced Mathematics. CRC Press, Boca Raton, FL}
\yr{1992}
\endbook

\bibitem{EH}
\by{\name{Esmayli}{B.} and \name{Hajlasz}{P.}}
\paper{The coarea inequality}
\jour{Ann. Fenn. Math.}
\vol{46}
\pages{965–991}
\yr{2021}
\endpaper



\bibitem{FG_95}
\by{\name{Fonseca}{I.} and \name{Gangbo}{W.}}
\paper{Local invertibility of Sobolev functions}
\jour{SIAM J. Math. Anal.}
\vol{26}
\yr{1995}
\pages{280--304}
\endpaper



\bibitem{HeMo11}
\by{\name{Henao}{D.} and \name{Mora-Corral}{C.}}
\paper{Fracture surfaces and the regularity of inverses for {BV} deformations}
\jour{Arch. Rational Mech. Anal.}
\vol{201}
\pages{575--629}
\yr{2011}
\endpaper
  
  \bibitem{HMC}
\by{\name{D.}{Henao} and \name{C.}{Mora-Corral}}
\paper{Lusin's condition and the distributional determinant for deformations with finite energy}
\jour{Adv. Calc. Var.}
\vol{5}
\pages{355--409}
\yr{2012}
\endpaper

  \bibitem{HMC2}
\by{\name{D.}{Henao} and \name{C.}{Mora-Corral}}
\paper{Regularity of inverses of Sobolev deformations with finite surface energy}
\jour{J. Funct. Anal.}
\vol{268}
\pages{2356--2378}
\yr{2015}
\endpaper

\bibitem{HK}
\by{\name{Hencl}{S.} and \name{Koskela}{P.}}
\book{Lectures on Mappings of finite distortion}
\publ{Lecture Notes in Mathematics 2096, Springer, 2014, 176pp}
\endbook

\bibitem{HK2}
\by{\name{Hencl}{S.} and \name{Koskela}{P.}}
\paper{Regularity of the inverse of a planar Sobolev homeomorphism}
\jour{Arch. Rational Mech. Anal}
\vol{180} 
\pages{75--95}
\yr{2006}
\endpaper


\bibitem{HKO05}
\by{\name{Hencl}{S.}, \name{Koskela}{P.} and \name{Onninen}{J.}}
\paper{A note on extremal mappings of finite distortion}
\jour{Math. Res. Lett.}
\vol{12}
\yr{2005}
\pages{231--237}
\endpaper


\bibitem{HKO}
\by{\name{Hencl}{S.}, \name{Koskela}{P.} and \name{Onninen}{J.}}
\paper{Homeomorphisms of bounded variation}
\jour{Arch. Rational Mech. Anal}
\vol{186} 
\pages{351--360}
\yr{2007}
\endpaper






\bibitem{IMb}
\by{\name{Iwaniec}{T.} and \name{Martin}{G.}}
\book{Geometric Function Theory and Non-linear Analysis. Oxford Mathematical Monographs}
\publ{Oxford University Press}
\yr{2001}
\endbook

\bibitem{IMO14}
\by{\name{Iwaniec}{T.}, \name{Martin}{G.} and \name{Onninen}{J.}}
\paper{On minimisers of Lp-mean distortion}
\jour{Comput. Methods Funct. Theory}
\vol{14}
\yr{2014}
\pages{399--416}
\endpaper


\bibitem{IMO}
\by{\name{Iwaniec}{T.}, \name{Martin}{G.} and \name{Onninen}{J.}}
\paper{Energy-minimal principles in geometric function theory}
\jour{New Zealand J. Math.}
\vol{52}
\yr{2021 [2021–2022]}
\pages{605--642}
\endpaper



\bibitem{IOmono}
\by{\name{Iwaniec}{T.} and \name{Onninen}{J.}}
\paper{Monotone Sobolev Mappings of Planar Domains and Surfaces}
\jour{Arch. Rational Mech. Anal.}
\vol{219}
\pages{159--181}
\yr{2016}
\endpaper



 \bibitem{K}
 \by{\name{Kalayanamit}{P.}}
 \paper{Sobolev regularity of the inverse for minimizers of the neo-Hookean energy satisfying condition $INV$}
 \jour{arXiv 2405.12156}
 \endprep





\bibitem{MM}
\by{\name{Marcus}{M.} and \name{Mizel}{J.}}
\paper{Transformations by functions in Sobolev spaces and lower semicontinuity for parametric variational problems}
\jour{Bull. Amer. Math. Soc.}
\vol {79, \rm no. 1}
\yr {1973}
\pages {790--795}
\endpaper

\bibitem{MY22}
\by{\name{Martin}{G. J.} and \name{Yao}{C.}}
\paper{The Teichmüller problem for $L^p$-means of distortion}
\jour{Ann. Fenn. Math.}
\vol{47}
\yr{2022}
\pages{1099--1108}
\endpaper


\bibitem{MY24}
\by{\name{Martin}{G.} and \name{Yao}{C.}}
\paper{The $L^p$ Teichm\"uller theory: existence and regularity of critical points}
\jour{Arch. Ration. Mech. Anal.}
\vol{248}
\yr{2024}
\pages{Paper No. 13, 35 pp}
\endpaper


\bibitem{Mor}
\by{\name{Morrey}{C. B.}}
\paper{The Topology of (Path) Surfaces}
\jour{Amer. J. Math.}
\vol{57}
\yr{1935}
\pages{17--50}
\endpaper


\bibitem{MS}
\by{\name{M\"uller}{S.} and \name{Spector}{S. J.}}
\paper{An existence theory for nonlinear elasticity that
allows for cavitation}
\jour{Arch. Ration. Mech. Anal.}
\vol {131, \rm no. 1}
\yr {1995}
\pages {1--66}
\endpaper

\bibitem{MST}
\by{\name{M\"uller}{S.}, \name{Spector}{S. J.} and \name{Tang}{Q.}}
\paper{Invertibility and a topological property of Sobolev maps}
\jour{SIAM J. Math. Anal.}
\vol{27}
\yr{1996}
\pages{959--976}
\endpaper

\bibitem{On06}
\by{\name{Onninen}{J.}}
\paper{Regularity of the inverse of spatial mappings with finite distortion}
\jour{Calc. Var. Partial Differential Equations}
\vol{26}
\yr{2006}
\pages{331--341}
\endpaper



\bibitem{Reb}
\by{\name{Reshetnyak}{Y. G.}}
\book{Space mappings with bounded distortion}
\publ{American Mathematical Society, Providence, RI}
\yr{1989}
\endbook



\bibitem{Sv}
\by{\name{\v Sver\'ak}{V.}}
\paper{Regularity properties of deformations with finite energy}
\jour{Arch. Ration. Mech. Anal.}
\vol{100\nom 2}
\pages{105--127}
\yr{1988}
\endpaper

\bibitem{SwaZie2002}
\by{\name{Swanson}{D.}, \name{Ziemer}{W. P.}}
\paper{A topological aspect of Sobolev mappings}
\jour{Calc. Var. Partial Differential Equations}
\vol{14\nom 1}
\pages{69--84}
\yr{2002}
\endpaper



\bibitem{SwaZie2004}
\by{\name{Swanson}{D.}, \name{Ziemer}{W. P.}}
\paper{The image of a weakly differentiable mapping}
\jour{SIAM J. Math. Anal.}
\vol{35\nom 5}
\pages{1099--1109}
\yr{2004}
\endpaper


\bibitem{T}
\by{\name{Tang}{Q.}}
\paper{Almost-everywhere injectivity in nonlinear elasticity}
\jour{Proc. Roy. Soc. Edinburgh Sect. A}
\vol{109}
\pages{79--95}
\yr{1988}
\endpaper



\bibitem{Vab}
J. V\"ais\"al\"a,  {\em Lectures on $n$-dimensional
quasiconformal mappings.} Lecture Notes in Mathematics, Vol. 229. 1971,
Springer-Verlag, Berlin-New York.


\bibitem{Z}
\by{\name{Ziemer}{W. P.}}
\book{Weakly differentiable functions. Graduate texts in Mathematics, 120}
\publ{Springer-Verlag}
\yr{1989}
\endbook




\end{thebibliography}
\end{document}